\documentclass[12pt, reqno]{amsart}
\usepackage{amsmath}
\usepackage{amssymb}
\usepackage{epsfig}
\usepackage{graphicx}
\usepackage{color}
\usepackage{fullpage}
\usepackage{comment}
\definecolor{shadecolor}{gray}{0.875}
\usepackage{amscd}
\usepackage{stmaryrd}

\numberwithin{equation}{section}

\input xy
\xyoption{all}

\calclayout
\allowdisplaybreaks[3]

\theoremstyle{plain}
\newtheorem{prop}{Proposition}[section]

\newtheorem{theo}[prop]{Theorem}
\newtheorem{coro}[prop]{Corollary}

\newtheorem{lemm}[prop]{Lemma}

\theoremstyle{definition}
\newtheorem{defi}[prop]{Definition}

\newtheorem{ques}[prop]{Question}
\newtheorem{conj}[prop]{Conjecture}

\newtheorem{rema}[prop]{Remark}

\newtheorem{exam}[prop]{Example}

\def\Eff{\overline{\mathrm{Eff}}}

\def\et{\mathrm{\acute{e}t}}

\def\Chow{\mathrm{Chow}}

\def\Nef{\mathrm{Nef}}

\def\Sym{\mathrm{Sym}}

\def\Sing{\mathrm{Sing}}
\def\Span{\mathrm{Span}}
\def\Gal{\mathrm{Gal}}
\def\Aut{\mathrm{Aut}}
\def\Jac{\mathrm{Jac}}

\makeatother
\makeatletter

\author{Brian Lehmann}
\address{Department of Mathematics \\
Boston College  \\
Chestnut Hill, MA \, \, 02467}
\email{lehmannb@bc.edu}

\author{Sho Tanimoto}
\address{Department of Mathematics, Faculty of Science, Kumamoto University, Kurokami 2-39-1 Kumamoto 860-8555 Japan}
\address{Priority Organization for Innovation and Excellence, Kumamoto University}
\email{stanimoto@kumamoto-u.ac.jp}

\newcounter{enumi_saved}

\title[Examples]{On exceptional sets in Manin's conjecture}

\begin{document}
\date{\today}

\begin{abstract}
In this survey paper we study Manin's Conjecture from a geometric perspective.  The focus of the paper is the recent conjectural description of the exceptional set in Manin's Conjecture due to Lehmann-Sengupta-Tanimoto.  After giving an extensive background, we give a precise description of this set and compute it in many examples. 
\end{abstract}

\maketitle

\section{Introduction}

Let $X$ be a smooth projective Fano variety over a number field $F$.  Let $\mathcal L = (L, \|\cdot\|)$ be an adelically metrized divisor on $X$. Then one can define a height function associated to $\mathcal L$
\[
\mathsf H_{\mathcal L} : X(F)\rightarrow \mathbb R_{\geq 0}.
\]
When the underlying divisor $L$ is ample, the height function satisfies the Northcott property, i.e., 
for any positive real number $T$, the set
\[
\{ P \in X(F) \mid \mathsf H_{\mathcal L} (P) \leq T \}
\]
is finite.  Thus, for any subset $Q \subset X(F)$, one can define the counting function
\[
N(Q, \mathcal L, T) := \#\{ P \in Q \mid \mathsf H_{\mathcal L} (P) \leq T \}.
\]
Manin's Conjecture, formulated in \cite{FMT89} and \cite{BM}, predicts that the asymptotic behavior of the counting function as $T$ increases is controlled by geometric invariants of $X$ and $L$.

Before stating the conjecture, we need to recall a definition: a thin subset of $X(F)$ is a finite union $\cup_{i} f_{i}(Y_{i}(F))$, where each $f_{i}: Y_{i} \to X$ is a morphism that is generically finite onto its image and admits no rational section.  Conceptually it is helpful to separate into two cases:
\begin{enumerate}
\item When $f_{i}$ is not dominant, then $f_{i}(Y_{i}(F))$ is not Zariski dense.
\item When $f_{i}$ is dominant, it is possible for $f_{i}(Y_{i}(F))$ to be Zariski dense, but it is still ``small'' in the same sense that the set of squares is ``small'' in the set of rational numbers.
\end{enumerate}
Thus we can view a thin set as a ``small'' set of points that generalizes the notion of a non-Zariski dense set of points.

Manin's Conjecture predicts that the growth of rational points is controlled by two geometric invariants $a(X,L)$ and $b(F,X,L)$ associated to $X$ and $L$.  We defer the definition of these invariants to Sections \ref{sect: asection} and \ref{sect: bsection}, where we discuss these two constants in great detail.

\begin{conj}[Manin's Conjecture] \label{conj: maninsconjecture}
Let $F$ be a number field.  Let $X$ be a geometrically rationally connected and geometrically integral smooth projective variety defined over $F$ and let $\mathcal{L}$ be a big and nef line bundle with an adelic metrization on $X$.

Suppose that $X(F)$ is not a thin set.
Then there exists a subset $Z \subset X(F)$ which is contained in a thin subset of $X(F)$ such that
\[
N(X(F) \setminus Z, \mathcal L, T) \sim c(F, Z, L)T^{a(X, L)} \log (T)^{b(F, X, L)-1}
\]
as $T \rightarrow \infty$ where $c(F, Z, L)$ is Peyre-Batyrev-Tschinkel's constant introduced in \cite{Peyre} and \cite{BT}.
\end{conj}

The subset $Z \subset X(F)$ in Manin's Conjecture is known as the exceptional set.  The original versions of Manin's Conjecture predicted that it suffices to remove an exceptional set of rational points which is not Zariski dense (\cite{BM} and \cite{Peyre}).  However, a counterexample to this version was identified in \cite{BT-cubic}; several more counterexamples to this version have been constructed since.  \cite{Peyre03} was the first to suggest that the exceptional set in Manin's Conjecture should be contained in a thin subset of rational points, and \cite{LeRudulier} was the first to exhibit such an instance of Manin's Conjecture.   There is currently no known counterexample to the thin set version of Manin's Conjecture.

In \cite{LST18}, the authors and A. Sengupta proposed a conjectural construction of the exceptional set $Z$.  This construction is motivated by geometry.  Suppose we have a morphism of smooth projective varieties $f: Y \to X$ which is generically finite and does not admit a rational section.  If
\begin{equation*}
(a(X,L),b(F,X,L)) < (a(Y,f^{*}L),b(F,Y,f^{*}L))
\end{equation*}
in the lexicographic order, then we expect the rate of growth of points on $Y$ to be higher than the desired growth rate on $X$.  Thus one should include $f(Y(F))$ in the exceptional set.  If we have equality of the $a,b$ invariants, then $f(Y(F))$ can affect the leading constant and one must consider carefully whether or not such points should be part of the exceptional set.  Loosely speaking, \cite{LST18} predicts that the exceptional set should be the union of the sets $f(Y(F))$ as we vary over all morphisms $f:Y \to X$ as above.  The main theorem of \cite{LST18} is that the proposed exceptional set is indeed contained in a thin set as predicted by Manin's Conjecture.

The goals of this paper are:
\begin{enumerate}
\item to give an introduction to the geometry of the $a$ and $b$ invariants.  In Sections \ref{sect: asection} and \ref{sect: bsection} we study their basic computational and theoretical properties.
\item to describe precisely the proposed exceptional set in \cite{LST18}.   We do this in Section \ref{sect: maninsconj}.  In particular, we give a careful accounting of when one should remove the point contributions from a morphism $f: Y \to X$ where the $a$ and $b$ invariants are equal.
\item to explain how to construct a thin set containing the exceptional set.  In Section \ref{sect: computation} we give an algorithm intended for arithmetic geometers who would like to compute the set constructed by \cite{LST18} in specific examples.  We also survey the theoretical results needed for proving thinness.
\item to illustrate the above results via examples, collected in Sections \ref{sect: hypersurfaces} -- \ref{sect: bhb}.
\end{enumerate}
Most of the techniques and theorems we describe are currently distributed in a series of papers \cite{HTT15}, \cite{LTT14}, \cite{HJ16}, \cite{LTDuke}, \cite{LT17}, \cite{Sen17b}, \cite{Sen17}, \cite{LST18}, and \cite{LT18}.  We present a few new results scattered throughout Sections \ref{sect: asection} and \ref{sect: bsection}.  Some of the examples in the later sections are also new.  We highlight three here:

\begin{exam}
In Section \ref{sect: surfaces} we prove that for a del Pezzo surface over a number field with the anticanonical polarization the proposed exceptional set of \cite{LST18} is not Zariski dense.  (When the del Pezzo surface has degree $1$ we need to assume that either the Picard rank is at least $2$ or the geometric model is general in moduli.)  This matches up with the currently known results on Manin's Conjecture for surfaces; see, e.g., \cite{dBBD07}, \cite{Bro09}, \cite{Bro10}, and \cite{dBBP12}.
\end{exam}

\begin{exam}
Suppose that $X$ is a toric variety and $L$ is a big and nef divisor on $X$.  We show that the only maps $f: Y \to X$ whose point contributions are removed by \cite{LST18} have image in the toric boundary.  Manin's Conjecture for toric varieties was proved in \cite{BT-general}, \cite{BT-0}, and \cite{Sal98} after removing rational points on the boundary.
\end{exam}

\begin{exam}
Recently \cite{BHB18} proved the thin set version of Manin's Conjecture for a diagonal quadric surface bundle.  In Section \ref{sect: bhb} we show that the proposed exceptional set of \cite{LST18}  exactly coincides with the thin exceptional set computed by \cite{BHB18} in this example.
\end{exam}

We close the paper with some open questions in Section \ref{sect: openques}.

\bigskip

\noindent
{\bf Acknowledgments.}
The authors would like to thank Marta Pieropan and Yuri Tschinkel for a stimulating question leading to this paper, and to thank Marta for many helpful comments on an earlier draft.
The authors would also like to thank Brendan Hassett, Akash Sengupta, and Yuri Tschinkel for collaborations helping to shape our perspective on the $a$ and $b$ invariants.
The first author would like to thank Jian Xiao for a useful conversation about \cite{KO73}.
The authors would like to thank the referee for careful reading of this paper and helpful suggestions.
Brian Lehmann is supported by NSF grant 1600875.  Sho Tanimoto is partially supported by MEXT Japan, Leading Initiative for Excellent Young Researchers (LEADER).

\section{Background}

Throughout we will work over a fixed field $F$ of characteristic $0$.  All morphisms will be $F$-morphisms unless otherwise noted.  We will denote an algebraic closure of $F$ by $\overline{F}$.  If $X$ is a scheme over $F$, we will denote by $\overline{X}$ its base change to $\overline{F}$.

A variety is a separated integral scheme of finite type over the ground field.

\begin{defi}
Let $X$ and $Y$ be projective varieties.  A map $f: Y \to X$ is thin if it is generically finite onto its image and admits no rational section.
\end{defi}

Thus if $F$ is a number field, a thin subset of $X(F)$ is precisely a finite union $\cup_{j} f_{j}(Y_{j}(F))$ where $f_{j}: Y_{j} \to X$ are thin maps over $F$.

\subsection{Pseudo-effective and nef cones}

On any projective variety $X$ there is an intersection pairing between Cartier divisors and $1$-cycles (that is, integer linear combinations of $1$-dimensional subvarieties of $X$).  We say that a Cartier divisor or a $1$-cycle is numerically trivial if it defines the zero linear functional under this pairing.  The quotient of the group of Cartier divisors by the subgroup of numerically trivial divisors is called the N\'eron-Severi group; we denote this group by $N^{1}(X)_{\mathbb{Z}}$.  This is a finite rank free abelian group.  The dual group of $1$-cycles up to numerical equivalence is denoted by $N_{1}(X)_{\mathbb{Z}}$.  We then set
\begin{equation*}
N^{1}(X) := N^{1}(X)_{\mathbb{Z}} \otimes \mathbb{R} \qquad \qquad N_{1}(X) := N_{1}(X)_{\mathbb{Z}} \otimes \mathbb{R}
\end{equation*}
The former is known as the N\'eron-Severi space of $X$.

Loosely speaking, one can view $N^{1}(X)$ as the natural ``homology space'' for algebraic divisors.  When our ground field is $\mathbb{C}$ and $X$ is smooth, then $N^{1}(X)$ is exactly the subspace of $H^{2}(X,\mathbb{R})$ spanned by divisor classes.  Analogously, $N_{1}(X)$ is the natural ``homology space'' for algebraic curves and for smooth projective varieties over $\mathbb{C}$ it coincides with the subspace of $H_{2}(X,\mathbb{R})$ spanned by curve classes.

\begin{exam}
If $X$ is a smooth projective Fano variety then numerical equivalence of Cartier divisors on $X$ coincides with linear equivalence and $N^{1}(X)_{\mathbb{Z}}$ is simply the Picard group of $X$.
\end{exam}

\begin{rema}
If $L$ is a Cartier divisor on $X$, we will also denote its equivalence class in $N^{1}(X)$ by $L$.  While a bit sloppy, in practice this convention is unlikely to cause confusion.
\end{rema}

The pseudo-effective cone of divisors, denoted by $\Eff^{1}(X)$, is the closure in $N^{1}(X)$ of the convex cone generated by classes of effective divisors.  Loosely speaking, one can view $\Eff^{1}(X)$ as the ``homological shadow'' of the algebraic structure of $X$ -- it records the possible homology classes of codimension $1$ subvarieties.  A divisor class is said to be big if it lies in the interior of $\Eff^{1}(X)$.  Analogously, the pseudo-effective cone of curves, denoted by $\Eff_{1}(X)$, is the closure in $N_{1}(X)$ of the convex cone generated by classes of effective curves.

The nef cone of divisors $\Nef^{1}(X)$ is the dual cone of the pseudo-effective cone of curves under the intersection pairing.  The famous result of \cite{kleiman66} shows that the nef cone of divisors is the same as the closure of the cone generated by classes of ample divisors.  Analogously, the nef cone of curves $\Nef_{1}(X)$ is the dual cone of the pseudo-effective cone of divisors.  By \cite{BDPP} it is the same as the closure of the cone generated by irreducible curves whose deformations dominate $X$.

\subsection{Big and nef divisors}

While the most important notion of ``positivity'' for divisors is ampleness, it will be helpful to have a slight generalization which is invariant under pullback by dominant generically finite morphisms.  A divisor class is said to be big and nef if it lies both in $\Nef^{1}(X)$ and in the interior of $\Eff^{1}(X)$.  Every ample class is big and nef, but the converse is not true.  As suggested above, the advantage of this notion is that the pullback of a big and nef divisor by a dominant generically finite morphism is again big and nef.

\begin{rema}
Sometimes it is useful to be a little more precise about the relationship between pullbacks and ample divisors, namely:
\begin{itemize}
\item The pullback of an ample divisor under a dominant generically finite map $f: X \to Y$ is a big and semiample divisor.  (The converse is also true if one works with $\mathbb{Q}$-divisors and $\mathbb{Q}$-linear equivalence.)
\item A divisor $L$ is big and basepoint free if and only if there is a dominant generically finite map $f: X \to Y$ and a very ample divisor $A$ on $Y$ such that $f^{*}A$ is linearly equivalent to $L$.
\end{itemize}
\end{rema}

An important distinction between ample divisors and big and nef divisors is that in the latter case it is possible for subvarieties to have vanishing $L$-degree.  However, the behavior of such subvarieties is controlled:

\begin{theo}[\cite{nakamaye00}] \label{theo: augbaselocus}
Let $X$ be a projective variety and let $L$ be a big and nef divisor on $X$.  The union of all subvarieties $V$ on $X$ satisfying $L^{\dim V} \cdot V = 0$ is a proper closed subset of $X$.
\end{theo}

When $X$ is smooth, the proper closed subset in Theorem \ref{theo: augbaselocus} is known as the augmented base locus of $L$ and is denoted by $\mathbf{B}_{+}(L)$.  Note that for an ample divisor $L$ the set $\mathbf{B}_{+}(L)$ is empty.  It turns out that many well-known theorems for ample divisors have analogues for big and nef divisors $L$.  For example, the set of subvarieties of $X$ of bounded $L$-degree which are not contained in $\mathbf{B}_{+}(L)$ is parametrized by a bounded subset of $\Chow(X)$.  Similarly, if we equip a big and nef line bundle $\mathcal{L} = \mathcal{O}_{X}(L)$ with an adelic metric, then the corresponding height function satisfies the Northcott property outside of $\mathbf{B}_{+}(L)$: for any positive integer $T$, there are only finitely many rational points in the complement of $\mathbf{B}_{+}(L)$ of $\mathcal{L}$-height $\leq T$.

\subsection{Canonical divisor}

Suppose that $X$ is a smooth variety of dimension $n$.  Then $\Omega_{X}$ denotes the sheaf of differential forms on $X$ which is a locally free sheaf of rank $n$.  Its top exterior power is called the canonical bundle and is denoted by $\omega_{X}$.  We will participate in the traditional abuse of notation by letting $K_{X}$ denote any Cartier divisor satisfying $\mathcal{O}_{X}(K_{X}) \simeq \omega_{X}$.  We call $K_{X}$ ``the'' canonical divisor on $X$.  Such a divisor is only unique up to linear equivalence, but since our statements are all linear-equivalence invariant this abuse is harmless in practice.

When $X$ is normal, one can still define the canonical divisor as a Weil divisor: there is a smooth subset $U \subset X$ whose complement has codimension $2$, and we define $K_{X}$ to be the unique Weil divisor corresponding to $K_{U}$.  However, when $X$ has worse singularities there may no longer be a natural way to define the canonical divisor.

There are various singularity types whose definitions involve the birational behavior of the canonical divisor: terminal, canonical, kawamata log terminal, etc.  The definitions can be found in any standard text on the minimal model program, e.g.~\cite{KM98}.

\section{The $a$-invariant} \label{sect: asection}

\subsection{Definition}
Let $X$ be a smooth projective variety and let $L$ be a big and nef divisor on $X$.  Our first invariant measures the negativity of $K_{X}$ against $L$.  Since $L$ is in the interior of $\Eff^{1}(X)$, we know that $K_{X} + mL$ will be contained in $\Eff^{1}(X)$ for all $m$ sufficiently large.  The $a$-invariant is the infimum of all such $m$:

\begin{defi}
Let $X$ be a smooth projective variety and let $L$ be a big and nef divisor.  We define the Fujita invariant, or the $a$-invariant, to be
\begin{equation} \label{eq: ainv}
a(X, L) = \min \{ t \in \mathbb R \mid K_{X} + tL \in \overline{\mathrm{Eff}}^1(X)\}.
\end{equation}
When $L$ is nef but not big, we formally set $a(X, L) = +\infty$.  (Note that we use this convention even when the right hand side of Equation \eqref{eq: ainv} is well-defined and finite.  This convention harmonizes with the fact that such $L$ can fail the Northcott property.)
\end{defi}

\begin{exam}
Suppose that $X$ is a Fano variety and that $L = -K_{X}$.  Then $a(X,L) = 1$.  Indeed, $K_{X} + tL$ will be anti-ample for $t<1$ and will be ample for $t>1$, and the boundary point is exactly $t=1$.
\end{exam}

\begin{exam}
Let $X$ be $\mathbb{P}^{n}$ and let $L$ be the hyperplane class.  Since $K_{X} = -(n+1)L$, we have $a(X,L) = n+1$.  
\end{exam}

\begin{exam}
Let $X$ be a smooth hypersurface of degree $d$ in $\mathbb{P}^{n}$ and let $L$ be the hyperplane class.  By adjunction we have $a(X,L) = n+1-d$.
\end{exam}

In examples where one can identify the pseudo-effective cone of divisors, it is usually easy to compute $a$-values by hand.

\begin{exam}
Let $X$ be $\mathbb{P}^{1} \times \mathbb{P}^{1}$ and let $H_{1}$ and $H_{2}$ denote the pullbacks of the hyperplane class under the first and second projections respectively.  With this notation we can identify $K_{X} = -2H_{1}-2H_{2}$.

In order to compute the $a$-invariant we must know the structure of the pseudo-effective cone.  In this situation $N^{1}(X)$ is spanned by $H_{1}$ and $H_{2}$ and
\begin{equation*}
\Eff^{1}(X) = \Nef^{1}(X) = \mathbb{R}_{\geq 0}H_{1} + \mathbb{R}_{\geq 0} H_{2}.
\end{equation*}
Any big and nef divisor $L$ will have the form $L = rH_{1} + sH_{2}$ where $r,s > 0$.  Then
\begin{equation*}
a(X,rH_{1} + sH_{2}) = \max \left\{ \frac{2}{r}, \frac{2}{s} \right\}.
\end{equation*}
\end{exam}

\begin{exam} \label{exam: blowupp2}
Let $X$ be the blow-up of $\mathbb{P}^{2}$ at a point.  Let $H$ denote the pullback of the hyperplane class on $\mathbb{P}^{2}$ and let $E$ denote the exceptional divisor.  With this notation we can identify $K_{X} = -3H + E$.

In this situation $N^{1}(X)$ is spanned by $H$ and $E$.  We have
\begin{align*}
\Eff^{1}(X) & = \mathbb{R}_{\geq 0}E + \mathbb{R}_{\geq 0}(H-E)  \\
\Nef^{1}(X) & = \mathbb{R}_{\geq 0}H + \mathbb{R}_{\geq 0} (H-E)
\end{align*}
Any big and nef divisor $L$ will have the form $L = rH + s(H-E)$ where $r > 0$ and $s \geq 0$.  Then
\begin{equation*}
a(X,(r+s)H - sE) = \max \left\{ \frac{3}{r+s}, \frac{2}{r} \right\}
\end{equation*}
\end{exam}

\subsection{First properties}
Note that the $a$-invariant is positive, zero, or negative when $K_{X}$ is contained in the complement, the boundary, or the interior of $\Eff^{1}(X)$ respectively.  We will mostly be interested in the first case, which is given a geometric interpretation by \cite{BDPP}.

\begin{theo}[\cite{BDPP}, 0.3 Corollary]
Let $X$ be a smooth projective variety and let $L$ be a big and nef divisor on $X$.  Then $a(X,L) > 0$ if and only if $X$ is geometrically uniruled.
\end{theo}

Using \cite{BCHM} we can reinterpret the $a$-constant using sections.  If $a(X,L)>0$ then we have an equality
\begin{equation} \label{eq: ainv2}
a(X, L) = \min \left\{ \left. \frac{r}{s} \in \mathbb Q \, \right| \dim H^{0}(X,sK_{X} + rL) > 0 \right\}.
\end{equation}
Note that this is a minimum, not an infimum; in particular $a(X,L)$ is always a rational number (but it need not be an integer).  The following list summarizes the key properties of the $a$-invariant.

\begin{prop}
Let $X$ be a smooth projective geometrically uniruled variety and let $L$ be a big and nef divisor on $X$.
\begin{itemize}
\item The $a$-invariant is $(-1)$-homogeneous: for any $t>0$ we have $a(X,tL) = \frac{1}{t} a(X,L)$.
\item The $a$-invariant is birational: if $\phi: X' \to X$ is a birational morphism of smooth projective varieties then $a(X,L) = a(X',\phi^{*}L)$.  (\cite[Proposition 2.7]{HTT15})
\item The $a$-invariant is geometric: if $F'/F$ is any field extension then $a(X,L) = a(X_{F'},L)$.
\end{itemize}
\end{prop}

\begin{rema}[Extension to singular varieties]
Suppose now that $X$ is a singular projective variety and that $L$ is a big and nef Cartier divisor on $X$.  It is not possible to define $a(X,L)$ using Equation \ref{eq: ainv} since for singular varieties it may not be possible to define $K_{X}$.  Even when $X$ is normal and $\mathbb{Q}$-Gorenstein, the presence of singularities will affect the computation in an undesirable way.

The simplest and best solution is to use  birational invariance: we define $a(X,L)$ by taking any resolution of singularities $\phi: X' \to X$ and setting $a(X,L) := a(X',\phi^{*}L)$.  This convention harmonizes nicely with the birational invariance of Manin's Conjecture.  When $X$ has only mild (i.e.~canonical $\mathbb{Q}$-Gorenstein) singularities, the right hand side of Equation \ref{eq: ainv} can still be used to compute $a(X,L)$.
\end{rema}

The following definition plays a crucial role in our analysis of the $a$-invariant.

\begin{defi}
Let $X$ be a smooth projective variety and let $L$ be a big and nef divisor on $X$.  We say that $(X,L)$ is adjoint rigid if for any sufficiently divisible positive integer $m$ we have that $\dim H^{0}(X,m(K_{X} + a(X,L)L)) = 1$.
\end{defi}

By an abuse of notation we will sometimes say that $(X,L)$ is adjoint rigid when $X$ is singular if the pullback of $L$ to a resolution of singularities yields an adjoint rigid pair.
Loosely speaking, one can think of an adjoint rigid pair as a birational version of a Fano variety.  More precisely:

\begin{lemm}[\cite{BCHM}]
Let $X$ be a smooth projective variety and let $L$ be a big and nef divisor on $X$ such that $(X,L)$ is adjoint rigid.  Let $E$ denote the unique effective divisor linearly equivalent to $m(K_{X} + a(X,L)L)$ for sufficiently divisible $m$.  Then there is a birational contraction map $\phi: X \dashrightarrow X'$ to a log Fano variety $X'$ which contracts $E$ to a subvariety of smaller dimension.
\end{lemm}

It follows from \cite{HM07} that if $X$ is a smooth geometrically integral projective variety, $L$ is a big and nef divisor, and $(X,L)$ is adjoint rigid then $X$ is geometrically rationally connected.

\subsection{Geometric properties of the $a$-invariant}

In this section we discuss how the $a$-invariant changes under some natural geometric operations.  The first and most important is the behavior of the $a$-invariant under covering maps.

\begin{lemm} \label{lemm: riemannhurwitz}
Let $X$ be a smooth projective variety and let $L$ be a big and nef divisor on $X$.  Then for any dominant generically finite morphism $f: Y \to X$ of smooth projective varieties we have $a(Y,f^{*}L) \leq a(X,L)$.
\end{lemm}

\begin{proof}
By the Riemann-Hurwitz formula there is some effective divisor $R$ such that $K_{Y} = f^{*}K_{X} + R$.  Since $K_{X} + a(X,L)L$ is pseudo-effective, we see that
\begin{equation*}
K_{Y} + a(X,L)f^{*}L = f^{*}(K_{X} + a(X,L)L) + R
\end{equation*}
is also pseudo-effective, proving the inequality of $a$-values.
\end{proof}

Using this result, one can prove:

\begin{lemm}[\cite{LTT14} Proposition 4.1]  \label{lemm: dominantfamilies}
Let $X$ be a smooth projective variety over an algebraically closed field carrying a big and nef divisor $L$ and let $Y$ be a general member of a dominant family of subvarieties of $X$.  Then $a(Y,L) \leq a(X,L)$.
\end{lemm}

Next, we discuss the behavior of the $a$-invariant in smooth families.  Using invariance of plurigenera as in \cite[Theorem 4.2]{HMX13} one can show:

\begin{theo} \label{theo: aconstantinfamilies}
Let $\pi: \mathcal{X} \to T$ be a smooth family of projective geometrically uniruled varieties.  Let $L$ be a big and semiample divisor on $\mathcal{X}$.  Then the function $T \to \mathbb{R}$ defined by $t \mapsto a(X_{t},L|_{X_{t}})$ is constant.
\end{theo}

When $L$ is only big and nef, \cite[Theorem 4.3]{LTDuke} proves that this function is lower semi-continuous, which is good enough for most applications.  (We are unsure whether one should expect the function to be constant in this situation.) Note that the theorem does not hold if we allow singular fibers in our family -- imagine for example a family of elliptic curves deforming to a nodal rational curve, where the $a$-invariant will be $0$ on the general fiber but positive on the singular fiber.
Finally, we discuss the behavior under taking hyperplane sections.

\begin{theo} \label{theo: hyperplanesection}
Let $X$ be a smooth projective geometrically uniruled variety of dimension $n$ and let $L$ be a big and basepoint free divisor on $X$.  Let $H$ be a general member of $|L|$.  Then one of the following holds:
\begin{enumerate}
\item $a(X,L) \leq 1$, or
\item $X$ is covered by geometrically rational curves $C$ with $L \cdot C = 1$, and we have $a(X,L)=2$ and $\kappa(K_{X} + a(X,L)L) = n-1$, or
\item $a(H,L) = a(X,L) - 1$ and $\kappa(K_{H} + a(H,L)L) = \kappa(K_{X} + a(X,L)L)$.
\end{enumerate}
\end{theo}

\begin{proof}  First we assume that the ground field is algebraically closed.
Suppose that $a(X,L) > 1$.  By \cite[Corollary D]{EP12} for $H$ general there is a surjection
\begin{equation*}
H^{0}(X,m(K_{X}+a(X,L)L)) \to H^{0}(H,m(K_{H}+(a(X,L)-1)L))
\end{equation*}
for any sufficiently divisible $m$. Furthermore the kernel of this map is $0$ since $K_{X} + a(X,L)L$ is on the pseudo-effective boundary and $H$ is big.

First suppose that $a(H,L) < a(X,L) - 1$ so that $K_{H} + (a(X,L)-1)L$ is big. Using the isomorphisms above we see that $\kappa(K_{X} + a(X,L)L) = n-1$.  Thus after making a birational change $\phi: X' \to X$ we obtain a morphism $\pi: X' \to Z$ whose general fiber is a rational curve $C'$ satisfying $a(C',\phi^{*}L) = a(X',\phi^{*}L)$. Let $C$ be $\phi_* C'$.  Since $a(C',\phi^{*}L) = 2/\phi^{*}L \cdot C'$, either $a(X,L) = a(C,L) \leq 1$ or $X$ is covered by rational curves with $L \cdot C = 1$ and $a(X,L)=a(C,L) = 2$.

Otherwise $a(H,L) = a(X,L) -1$ and we have an equality $\kappa(K_{H} + a(H,L)L) = \kappa(K_{X} + a(X,L)L)$.

If the ground field is not algebraically closed, then the only thing to check is the existence of the curves $C$ over the ground field in case (2).  But this follows from the fact that the map defined by sections of multiples of $K_{X} + a(X,L)L$ is defined over the ground field.
\end{proof}

\subsection{Varieties with large $a$-invariant}

The Fujita invariant was introduced by \cite{Som86} and \cite{Fuj92} due to its relationship with the positivity of divisors of the form $K_{X} + mA$ for an ample divisor $A$.  When $\dim H^{0}(X,K_{X} + mL) > 0$ we know that $m \geq a(X,L)$.  Note however that when $K_{X} + mL$ fails to have sections we can not deduce anything about $a(X,L)$ (see Equation \eqref{eq: ainv2}).  Thus the study of sections of $K_{X} + mL$ is mostly useful for constructing upper bounds on $a(X,L)$.  In this subsection we discuss some applications to varieties with large $a$-value.

\begin{lemm} \label{lemm: siu}
Let $X$ be a smooth projective variety of dimension $n$ and let $L$ be a big and nef divisor on $X$.  Then $a(X,L) \leq n+1$.
\end{lemm}

To the best of our knowledge this argument is due to Siu.

\begin{proof}
Set $P(m) := \chi(\mathcal{O}_{X}(K_{X} + mL))$.  By Hirzebruch-Riemann-Roch, 
$P(m)$ is a polynomial in $m$ of degree at most $n$.  By Kawamata-Viehweg vanishing, for $m>0$ we have
\begin{equation*}
P(m) = H^{0}(X,K_{X} + mL).
\end{equation*}
Note that $P(m)$ is not identically zero since for $m$ large  we have 
$H^{0}(X,K_{X}+mL) > 0$.  Thus $P(m)$ cannot have $n+1$ roots, so that $H^{0}(X,K_{X}+mL) > 0$ for some $1 \leq m \leq n+1$.
\end{proof}

The papers \cite{Fujita89}, \cite{Horing10}, \cite{Andreatta13} give a classification of varieties with large $a$-value in the spirit of the work of Kobayashi-Ochiai (\cite{KO73}).  The first step in the classification is:

\begin{prop}[\cite{Fujita89}, \cite{Horing10}] \label{prop: fujitaclassification}
Let $X$ be a smooth projective variety of dimension $n$ over an algebraically closed field and let $L$ be a big and nef divisor on $X$.
\begin{itemize}
\item If $a(X,L) > n$ then $a(X,L)=n+1$ and the pair $(X,L)$ is birationally equivalent to $(\mathbb{P}^{n},H)$ where $H$ is the hyperplane divisor on $\mathbb{P}^{n}$.
\item If $a(X,L) = n$ and $(X,L)$ is adjoint rigid then $(X,L)$ is birationally equivalent to $(Q,H)$ where $Q$ is a (possibly singular) quadric hypersurface and $H$ is the hyperplane class on $Q$.
\item If $a(X,L) = n$ and $(X,L)$ is not adjoint rigid, up to birational equivalence the canonical map $\pi: X \to C$ realizes $X$ as a $\mathbb{P}^{n-1}$-bundle over a curve $C$ and $L = \mathcal{O}_{\pi}(1)$.
\item If $n-1 < a(X,L) < n$ then $(X,L)$ is birationally equivalent to $\mathbb{P}_{\mathbb{P}^{2}}(\mathcal{O}(2) \oplus \mathcal{O}^{\oplus (n-2)})$ with the divisor $\mathcal{O}_{X/\mathbb{P}^{2}}(1)$.  In this case $(X,L)$ is adjoint rigid and $a(X,L) = n - \frac{1}{2}$.
\end{itemize}
\end{prop}

Here we say that two pairs $(X,L)$ and $(X',L')$ are birationally equivalent if there is a projective variety $W$ and birational maps $\phi_{1}: W \to X$ and $\phi_{2}: W \to X'$ such that $\phi_{1}^{*}L$ and $\phi_{2}^{*}L'$ are numerically equivalent.

More generally, one can use techniques from the study of Fujita's Conjecture to analyze $a$-values.  \cite{LTT14} gives some applications of Reider's Theorem to the study of $a$-values for surfaces.  It turns out that techniques arising from the minimal model program are usually more powerful, and we discuss these next.

\subsection{The $a$-invariant and the minimal model program}
\label{subsec:MMP}

Let $X$ be a smooth projective variety and let $L$ be a big and nef divisor on $X$.  The minimal model program associates to $K_{X} + a(X,L)L$ a rational map $\pi: X \dashrightarrow Z$ which encapsulates vital information about the $a$-invariant. 

\begin{theo}[\cite{BCHM} Theorem 1.1]
Let $X$ be a smooth geometrically uniruled projective variety and let $L$ be a big and nef divisor on $X$.  Choose a positive integer $d$ such that $d(K_{X} + a(X,L)L)$ is a Cartier divisor.  Then the section ring
\begin{equation*}
R(X,K_{X} + a(X,L)L) := \oplus_{m \geq 0} H^{0}(X,md(K_{X} + a(X,L)L))
\end{equation*}
is finitely generated.
\end{theo}

Since this ring is finitely generated, we can define a projective variety $Z = \mathrm{Proj} \, R(X,K_{X} + a(X,L)L)$.  This variety is known as the ``canonical model'' for the pair $(X,a(X,L)L)$.  Using the properties of the Proj construction we obtain a rational map $\pi: X \dashrightarrow Z$ known as the ``canonical map'' for the pair $(X,a(X,L)L)$.  (This is also the rational map defined by sections of any sufficiently large multiple of $K_{X} + a(X,L)L$.)  Note that $(X,L)$ is adjoint rigid if and only if $\pi$ contracts $X$ to a point.

In fact \cite{BCHM} tells us much more; we know that the rational map $\pi: X \dashrightarrow Z$ can be resolved by a sequence of birational changes known as ``running the minimal model program.''

\begin{lemm}[\cite{KM98} Lemma 3.38 and Corollary 3.53]
Let $X$ be a smooth uniruled projective variety and let $L$ be a big and nef divisor on $X$.  Then the $a$-invariant is preserved by birational steps of the $(X,a(X,L)L)$-MMP.  More precisely, if $\phi: X \dashrightarrow X'$ is a birational contraction constructed from steps of the minimal model program, then $a(X,L) = a(X',\phi_{*}L)$.
\end{lemm}

\begin{rema}
In the setting of the lemma above $\phi_{*}L$ need not be nef, but we nevertheless define $a(X',\phi_{*}L)$ using Equation \eqref{eq: ainv}.
\end{rema}

The most important consequences for us are summarized by the next proposition:

\begin{prop} \label{prop: coveredbyadjrigid}
Let $X$ be a smooth uniruled projective variety and let $L$ be a big and nef divisor on $X$.  Let $F$ be (the closure of) a general fiber of the canonical model map $\pi: X \dashrightarrow Z$ for $(X,a(X,L)L)$.  Then:
\begin{itemize}
\item We have $a(X,L) = a(F,L)$.
\item The pair $(F,L)$ is adjoint rigid.
\end{itemize}
\end{prop}

In other words, every variety $X$ with positive $a$-value is birationally fibered by adjoint rigid varieties with the same $a$-value.  In order to study the original variety, we can think of it as a family of adjoint rigid varieties.  We will use this idea pervasively in the following sections to reduce questions about arbitrary varieties to the adjoint rigid case.

\subsection{Boundedness results}

The main advantage of working with adjoint rigid pairs is that they have very special geometric properties.

\begin{prop}[\cite{LTDuke} Theorem 3.5]
Let $X$ be a smooth uniruled projective variety of dimension $n$ over an algebraically closed field and let $L$ be a big and nef divisor such that $(X,L)$ is adjoint rigid.  Then $X$ is birational to a $\mathbb{Q}$-Gorenstein Fano variety $X'$ with canonical singularities satisfying $a(X,L)^{n}L^{n} \leq (-K_{X'})^{n}$.
\end{prop}

This proposition allows us to leverage the recent solution of the Borisov-Alexeev-Borisov Conjecture in \cite{birkar16}, \cite{birkar16b}.  In these papers Birkar proves that all Fano varieties with mild singularities lie in a bounded family, and in particular, admit a universal upper bound for the anticanonical volume.  Thus:

\begin{theo} \label{theo: babbound}
There is a constant $C = C(n)$ satisfying the following property.  Suppose that $X$ is a smooth uniruled projective variety of dimension $n$ over an algebraically closed field and that $L$ is a big and nef divisor on $X$ such that $(X,L)$ is adjoint rigid.  Then $a(X,L)^{n}L^{n} \leq C$.
\end{theo}

While this theorem is not effective in higher dimensions, for surfaces and threefolds there are nearly-tight bounds on the anticanonical volume of mildly singular Fanos.  It is worth emphasizing that the adjoint rigid assumption is necessary -- there is in general no boundedness theorem for the $a$-invariant.

\begin{exam}
Let $X$ be the product of $\mathbb{P}^{1}$ with a non-uniruled surface $S$.  Let $L$ be any ample divisor on $X$.  If $C$ denotes a fiber of the projection $\pi: X \to S$, then $a(X,L) = 2/L \cdot C$ and $\pi$ will be the map to the canonical model for $(X,a(X,L)L)$.

Let $Y \subset X$ denote the $\pi$-preimage of any curve in $S$.  Then $a(Y,L) = a(X,L)$.  In particular, subvarieties of $X$ with the same $a$-value as $X$ can have arbitrarily large degree against $L$.
\end{exam}

\subsection{Computing the $a$-invariant in small dimensions}

For small dimensions we can apply previous results to give many concrete techniques for computing the $a$-invariant.

\subsubsection{Curves}

Suppose that $C$ is an irreducible (but not necessarily smooth) rational curve.  For any big and nef Cartier divisor $L$ on $C$ we have $a(C,L) = 2/\mathrm{deg}_{C}(L)$.

\subsubsection{Surfaces}
If we interpret Theorem \ref{theo: babbound} in this situation, we obtain the following explicit bound:

\begin{theo} \label{theo: surfacebabbound}
Let $S$ be a projective surface and let $L$ be a big and nef Cartier divisor on $X$.  Then either $S$ is dominated by curves $C$ satisfying $a(C,L) = a(S,L)$ or
\begin{equation*}
a(S,L)^{2}L^{2} \leq 9.
\end{equation*}
\end{theo}

The minimal model program can be used to derive more specific information.  We start from a smooth projective surface $S$ and a big and nef divisor $L$.  Suppose that $K_{S} + a(S,L)L$ is not nef.  Then there must be a Galois orbit of $(-1)$-curves $E$ on $S$ satisfying $a(S,L)L \cdot E < 1$.  Let $\phi: S \to S'$ be the contraction of these curves.  Then $L' = \phi_{*}L$ is still big and nef and we have $a(S,L) = a(S',L')$.

Continuing in this way we eventually find a birational model $\widetilde{\phi}: S \to \widetilde{S}$ for which $K_{\widetilde{S}} + a(S,L)\widetilde{L}$ is nef where $\widetilde{L} = \widetilde{\phi}_{*}L$.  There are two possibilities:
\begin{itemize}
\item Some multiple of this divisor defines a map $\pi: \widetilde{S} \to Z$ to a curve $Z$.  (This will of course be the canonical map discussed in Section~\ref{subsec:MMP}.)  In this case the general fiber of $\pi$ will be a curve with the same $a$-value as $S$.
\item We have $K_{\widetilde{S}} = -a(S,L)\widetilde{L}$.  In this case $\widetilde{S}$ is known as a weak del Pezzo surface.  This class of surfaces is well-understood, and the $a$-invariant can be computed easily once the surface $\widetilde{S}$ is identified and $\widetilde{L}^{2}$ is calculated.
\end{itemize}

Here are a couple theorems that can be proved using these techniques:

\begin{lemm}[\cite{LT18} Lemma 5.3] \label{lemm: normalsurfacebound}
Suppose that $Y$ is a normal surface and $L$ is an ample Cartier divisor on $Y$.  Then one of the following holds:
\begin{itemize}
\item $a(Y,L)<1$, or
\item $Y$ has canonical singularities and $a(Y,L)L \equiv -K_{Y}$, or
\item for a resolution $\phi: \widetilde{Y} \to Y$ the adjoint pair $K_{\widetilde{Y}} + a(Y,L)\phi^{*}L$ has Iitaka dimension $1$.
\end{itemize}
\end{lemm}

\begin{lemm}[\cite{LTDuke} Theorem 3.7 and \cite{LT18} Lemma 5.4]
Suppose that $Y$ is a normal surface and $L$ is an ample Cartier divisor on $Y$.
\begin{itemize}
\item Suppose that the smallest $L$-degree of a rational curve through a general point of $Y$ is $> 3$.  Then $a(Y,L) < 1$.
\item Suppose that the smallest $L$-degree of a rational curve through a general point of $Y$ is $3$.  Then either $a(Y,L) < 1$ or $Y \cong \mathbb{P}^{2}$ and $L = \mathcal{O}(1)$.
\end{itemize}
\end{lemm}

\subsubsection{Threefolds}  The cleanest statement for threefolds is:

\begin{theo}[\cite{LTDuke} Theorem 3.7]
Let $X$ be a smooth projective threefold and let $L$ be a big and nef divisor on $X$.  Then either $a(X,L)^{3} \leq 64/L^{3}$ or $X$ is dominated by subvarieties with the same $a$-value.
\end{theo}

Here the constant $64$ is arising from an explicit anticanonical volume bound for singular Fanos.  It is helpful to know the exact bound in various situations:

\begin{theo}
\begin{enumerate}
\item If $X$ is a Gorenstein terminal Fano threefold then $-K_{X}^{3} \leq 64$, and the bound is sharp. (\cite{Nami})
\item If $X$ is a $\mathbb Q$-factorial non-Gorenstein terminal Fano threefold of Picard rank $1$ then $-K_{X}^{3} \leq 125/2$, and the bound is sharp. (\cite{Pro})
\item If $X$ is a Gorenstein canonical Fano threefold then $-K_{X}^{3} \leq 72$, and the bound is sharp. (\cite{Pro2})
\end{enumerate}
\end{theo}

\begin{rema}
In dimension $\geq 4$ there are smooth Fano varieties whose anticanonical volume is larger than the anticanonical volume $(n+1)^{n}$ of $\mathbb{P}^{n}$.  However, projective space does have the largest anticanonical volume amongst all smooth K\"ahler-Einstein Fano varieties; see \cite{fujita15}. 
\end{rema}

\section{The $b$-invariant} \label{sect: bsection}

\begin{defi} \label{defi: binvariant}
Let $X$ be a smooth geometrically integral projective variety over a field $F$ and let $L$ be a big and nef divisor on $X$.  We define the $b$-invariant by
\begin{align*}
b(F, X, L) = \textrm{the }& \text{codimension of the minimal supported face}\\
& \text{of $\overline{\mathrm{Eff}}^1(X)$ containing $K_{X} + a(X, L)L$.}
\end{align*}
We will only ever consider this invariant when $a(X,L)>0$.  If $L$ is nef but not big, we formally set $b(X,L) = \infty$.
\end{defi}

Recall that a face $\mathcal{F}$ of a cone $\mathcal{C}$ is said to be supported if there is a linear functional $\ell$ such that $\ell(v) \geq 0$ for every $v \in \mathcal{C}$ and $\mathcal{F} = \{ v \in \mathcal{C} | \ell(v) = 0\}$. In general the issue of whether a face is supported or not can be somewhat subtle.  However, the situation is better for the face $\mathcal{F}$ of Definition \ref{defi: binvariant}.  According to \cite[Theorem 2.16]{HTT15} the minimal supported face of $\Eff^{1}(X)$ containing $K_{X} + a(X,L)L$ coincides with the minimal face of $\Eff^{1}(X)$ containing $K_{X} + a(X,L)L$.  Thus $K_{X} + a(X,L)L$ lies in the interior of the face.

\begin{exam} \label{exam: fanobinv}
Suppose that $X$ is a smooth Fano variety and that $L = -K_{X}$.  Then $b(F, X,L)$ is the Picard rank of $X$.
\end{exam}

\begin{exam}
We return to the setting of Example \ref{exam: blowupp2}.   Let $X$ be the blow-up of $\mathbb{P}^{2}$ at a point.  Let $H$ denote the pullback of the hyperplane class on $\mathbb{P}^{2}$ and let $E$ denote the exceptional divisor.  With this notation we can identify $K_{X} = -3H + E$.

As discussed in Example \ref{exam: blowupp2}, we have
\begin{equation*}
a(X,(r+s)H - sE) = \max \left\{ \frac{3}{r+s}, \frac{2}{r} \right\}.
\end{equation*}
It is straightforward to see that
\begin{equation*}
b(F,X,(r+s)H - sE) = \left\{ \begin{array}{ll} 2 & \textrm{ if }L \textrm{ is proportional to } -K_{X} \textrm{, i.e. }r=2s \\
1 & \textrm{ otherwise}\end{array} \right.
\end{equation*}
\end{exam}

\subsection{First properties}
We discuss some of the elementary properties of the $b$-invariant.

\begin{prop}
Let $X$ be a smooth geometrically integral geometrically uniruled projective variety and let $L$ be a big and nef divisor on $X$.
\begin{itemize}
\item The $b$-invariant is invariant under rescaling: for any $t>0$ we have $b(F,X,tL) = b(F,X,L)$.
\item The $b$-invariant is birational: if $\phi: X' \to X$ is a birational morphism of smooth projective varieties then $b(F,X,L) = b(F,X',\phi^{*}L)$.  (\cite[Proposition 2.10]{HTT15})
\end{itemize}
\end{prop}

However, as suggested by the notation the $b$-invariant is \emph{not} geometric: it can increase under base change (e.g.~when $X$ is Fano and the the Picard rank increases upon base change).

\begin{prop}
Let $X$ be a smooth geometrically integral geometrically uniruled projective variety defined over a field $F$ and let $L$ be a big and nef divisor on $X$. Let $F'/F$ be a finite extension. Then we have
\[
b(F, X, L) \leq b(F', X_{F'}, L_{F'}).
\]
\end{prop}

\begin{proof}
Let $\overline{F}$ be an algebraic closure of $F$.  Then $b(F, X, L)$ is given by
\[
\dim \overline{\mathrm{Eff}}_1(X_{\overline{F}})^{\mathrm{Gal}(\overline{F}/F)} \cap (K_{X} + a(X,L)L)^{\perp}.
\]
It is then clear that the $b$-value can only increase upon field extension.
\end{proof}

\subsection{The $b$-invariant and the minimal model program}

It is clear that the $b$-invariant is somehow related to the Picard rank of $X$ and the Picard rank of $Z$ where $\pi: X \dashrightarrow Z$ is the canonical map associated to $(X,a(X,L)L)$.  In this section we make this intuition precise.

It is most convenient to separate into two cases: when $(X,L)$ is adjoint rigid and when it is not.

\begin{lemm} \label{lemm: binvadjrigid}
Let $X$ be a smooth geometrically uniruled geometrically integral projective variety defined over $F$ and let $L$ be a big and nef divisor on $X$.  Suppose that $(X,L)$ is adjoint rigid.  Fix a sufficiently divisible integer $d$.  Then there is a unique effective divisor $E$ such that $md(K_{X} + a(X,L)L) \sim mE$ for all $m \geq 0$.  We have
\begin{equation*}
b(F,X,L) = \dim N^{1}(X)/\Span\{E_{i}\}
\end{equation*}
where $E_{i}$ are the irreducible components of $E$.
\end{lemm}

Note that $E$ is exactly the divisor that gets contracted when one runs the minimal model program on $X$.  If we denote the resulting birational map by $X \dashrightarrow X'$, then the $b$-value of $(X,L)$ is precisely the Picard rank of $X'$.

When $X$ is not adjoint rigid, the $b$-invariant is measured by a ``relative'' analogue of the calculation in Lemma \ref{lemm: binvadjrigid}.  In essence, the $b$-invariant measures the monodromy-invariant Picard rank of a general fiber of the canonical map after quotienting by the horizontal rigid divisors.

\begin{lemm}[\cite{LST18} Lemma 2.6] \label{lemm: monodromyandbvalue}
Let $X$ be a smooth geometrically uniruled geometrically integral projective variety defined over $F$. Let $\pi: X \dashrightarrow Z$ be the canonical map for $K_{X} + a(X,L)L$.  Suppose that there is a non-empty open set $Z^{\circ} \subset Z$ such that $\pi$ is well-defined and smooth over $Z^{\circ}$. For any $z \in Z^{\circ}(F)$ we have
\begin{align*}
b(F,X,L) & = \dim N^{1}(\overline{X}_{z})^{\pi^{\et}_1(Z^{\circ},z)} / \Span \{ E_{i} \}
\end{align*}
where the $\{ E_{i} \}$ are the classes of all irreducible divisors which dominate $Z$ and which satisfy $K_{X} + a(X,L)L - c_iE_{i} \in \Eff^{1}(X)$ for some $c_i > 0$.  Furthermore, there are only finitely many such $E_{i}$.
\end{lemm}


Another useful observation is that the $b$-invariant is preserved by the steps of the MMP:

\begin{lemm}[\cite{HTT15} Proposition 2.18] \label{lemm: binvandmmp}
Let $X$ be a smooth uniruled projective variety and let $L$ be a big and nef divisor on $X$.  Then the $b$-invariant is preserved by birational steps of the $(X,a(X,L)L)$-MMP.  More precisely, if $\phi: X \dashrightarrow X'$ is a birational contraction constructed from steps of the minimal model program, then $b(F,X,L) = b(F,X',\phi_{*}L)$.
\end{lemm}

\begin{rema}
In the setting of the lemma above $\phi_{*}L$ need not be nef, but we nevertheless define $b(F,X',\phi_{*}L)$ as in Definition \ref{defi: binvariant}.
\end{rema}

\begin{rema}
It is tempting to try to combine Lemma \ref{lemm: monodromyandbvalue} and Lemma \ref{lemm: binvandmmp} to measure the $b$-invariant using the fibers of the canonical map after running the MMP.  However, the analogue of Lemma \ref{lemm: monodromyandbvalue} may not be valid in this case if the general fiber of the canonical map fails to be $\mathbb{Q}$-factorial.  See \cite[Section 2]{CFST16} for an in-depth discussion of this issue.

\end{rema}

\subsection{Geometric properties of the $b$-invariant}

In this section we discuss how the $b$-invariant changes under some natural geometric operations.  Unfortunately it is not nearly as well behaved as the $a$-invariant.  The following examples show that the $b$-invariant can either increase or decrease under taking covers with the same $a$-invariant.

\begin{exam}
Suppose $F$ is algebraically closed.  Then \cite[Corollary 6.6]{LTDuke} shows that the $b$-invariant will always decrease for a covering map of adjoint rigid surfaces with the same $a$-invariant.  One can construct such maps by taking del Pezzo surfaces with canonical singularities and constructing covers which are \'etale in codimension $1$; some explicit examples can be found in \cite{MZ88}.
\end{exam}

\begin{exam}[\cite{LeRudulier}]
Let $S$ be the quadric surface $\mathbb P^1 \times \mathbb P^1$ defined over $\mathbb Q$ and consider the Hilbert scheme $X = \mathrm{Hilb}^{[2]}(S)$; this is a weak Fano $4$-fold. Let $L = -K_X$. Then we have
\[
a(X, L) = 1, \quad b(\mathbb Q, X, L) =  3,
\]
because the Picard rank of $X$ is $3$. Consider the variety $W$ which is the blow-up of $S \times S$ along the diagonal.  
These fit into a diagram
\begin{equation*}
\xymatrix{W \ar[r]\ar[d]_{f}&  S \times S \ar[d]^{g}\\
X \ar[r] & \Chow^{2}(S)}
\end{equation*}
where the vertical maps are the natural symmetric quotients.  Then we have 
\[
a(W,f^{*}L)=1 \quad b(\mathbb Q, W,f^{*}L) = 4
\]
and $(W, f^*L)$ is adjoint rigid.  Thus the $b$-value can increase for adjoint rigid covers with the same $a$-invariant.  
\end{exam}

Next, we discuss the behavior of the $b$-invariant in smooth families.

\begin{theo}[\cite{Sen17b} Theorem 1.2] \label{theo: bconstantinfamilies}
Let $\pi: \mathcal{X} \to T$ be a smooth family of projective uniruled varieties over an algebraically closed field $F$.  Let $L$ be a big and nef divisor on $\mathcal{X}$.  Consider the function $T \to \mathbb{R}$ defined by $t \mapsto b(F,X_{t},L|_{X_{t}})$.  There is an open subset $T^{\circ} \subset T$ such that the $b$-invariant is constant on $T^{\circ}$.
\end{theo}

When the ground field is not algebraically closed, the behavior can be more complicated.  For example, in a smooth family of Fano varieties over a number field the Picard rank can jump up over a thin dense set of points.  Thus if we take $L$ to be the relative anticanonical divisor, the $b$-value may not be constant on any open subset.

Finally, we discuss the behavior under taking hyperplane sections.

\begin{theo}
Let $X$ be a smooth projective uniruled variety of dimension $\geq 4$ and let $L$ be a very ample divisor on $X$.  Let $H$ be a general member of $|L|$.  Suppose that $a(X,L) > 1$ and that $a(H,L) = a(X,L)-1$.  Then $b(F,X,L) \geq b(F,H,L)$.

Suppose furthermore that there is an open neighborhood $U$ of $L$ contained in the ample cone such that $\kappa(X,K_{X} + a(X,L')L') < \dim(X)-1$ for any $\mathbb{Q}$-divisor $L'$ in $U$.  Then $b(F,X,L) = b(F,H,L)$.
\end{theo}

\begin{proof}
By the Lefschetz theorem the restriction map $\mathrm{res}: N^{1}(X) \to N^{1}(H)$ is bijective.  Furthermore $\mathrm{res}(\Eff^{1}(X)) \subset \Eff^{1}(H)$.  Since $K_{H} + (a(H,L))L = \mathrm{res}(K_{X} + a(X,L)L)$ we see that the minimal face containing the adjoint divisor can only increase in dimension.  Thus $b(F,X,L) \geq b(F,H,L)$.  

To show the second statement, note that after possibly shrinking $U$ we may suppose that for any $\mathbb{Q}$-divisor $L'$ in $U$ there is an ample divisor $A_{L'}$ such that
\begin{equation*}
a(X,L')L' = H + A_{L'}.
\end{equation*}
Using \cite[Corollary D]{EP12} and the assumption on Iitaka dimensions, by arguing as in Theorem \ref{theo: hyperplanesection} we see that $a(X,L') - 1 = a(H,L')$ for any $\mathbb{Q}$-divisor $L'$ in $U$.  In other words, there is an open neighborhood $V$ of $K_{X} + a(X,L)L$ such that
\begin{equation*}
\mathrm{res}(V \cap \Eff^{1}(X)) \cap N^{1}(H)_{\mathbb{Q}} = \mathrm{res}(V) \cap \Eff^{1}(H) \cap N^{1}(H)_{\mathbb{Q}}.
\end{equation*}
Since the minimal supported face of $\Eff^{1}(X)$ containing $K_{X} + a(X,L)L$ is rational and contains $K_{X} + a(X,L)L$ in its interior, this suffices to show equality of $b$-invariants.
\end{proof}

\section{The exceptional set in Manin's Conjecture} \label{sect: maninsconj}

Let $X$ be a smooth projective geometrically rationally connected and geometrically integral variety over a number field $F$ and let $L$ be a big and nef divisor on $X$.  Recall that a map of projective varieties $f: Y \to X$ is thin when $f$ is generically finite onto its image but is not birational and dominant.  \cite{LST18} predicts that the exceptional set in Manin's Conjecture is a union of subsets $f(Y(F))$ as $f: Y \to X$ varies over thin maps satisfying certain geometric conditions.  The goal of this section is to give a precise description of the set constructed in \cite{LST18}.

\subsection{Fano varieties}
The most important case is when $X$ is a smooth projective geometrically integral Fano variety over a number field and $L = -K_{X}$.  Suppose that $f: Y \to X$ is a thin map where $Y$ is a smooth projective geometrically integral variety.  As discussed in the introduction, the exceptional set in Manin's Conjecture should certainly include $f(Y(F))$ when we have an inequality
\[
(a(Y, f^*L), b(F, Y, f^*L)) > (a(X, L), b(F, X, L))
\]
in the lexicographic order.  If $f$ yields an equality then sometimes the contributions from $Y$ should be removed and sometimes they should not, and it is important to carefully distinguish between the two situations.

First suppose that $\dim(Y) < \dim(X)$ and the $a$ and $b$ invariants are equal.  Then we should always remove point contributions from $Y$.  Indeed, \cite[Theorem 1.2]{BL16} shows that if rational points grow at the expected rate on $Y$ and we include these points in our count then Manin's Conjecture with Peyre's constant will be violated for an appropriate choice of anticanonical height function.  In fact, we should remove the points on any family of such subvarieties, or even a cover of such a family.  Precisely: suppose that $f: Y \to X$ is a thin map satisfying
\[
(a(Y, f^*L), b(F, Y, f^*L)) = (a(X, L), b(F, X, L))
\]
and that the fibers of the canonical map for $(Y,f^{*}L)$ have smaller dimension than $X$.  Then we include $f(Y(F))$ in the set we must remove.

The other possibility is that $f: Y \to X$ is a dominant thin map satisfying 
\[
(a(Y, f^*L), b(F, Y, f^*L)) = (a(X, L), b(F, X, L))
\]
and $(Y,f^{*}L)$ is adjoint rigid.  In this situation, \cite{LT17} introduced a new invariant to distinguish when $f(Y(F))$ should be included.

\begin{defi}{\cite[Definition 3.5]{LT17}} \label{defi: facecontraction}
Let $X$ be a smooth geometrically integral projective variety and let $L$ be a big and nef divisor on $X$.  Suppose that $f: Y \to X$ is a morphism of smooth projective varieties that is dominant generically finite and satisfies $a(Y,f^{*}L) = a(X,L)$.  Let $F_{Y}$ denote the face of $\Nef_{1}(Y)$ perpendicular to $K_{Y} + a(Y,f^{*}L)f^{*}L$ and let $F_{X}$ denote the face of $\Nef_{1}(X)$ perpendicular to $K_{X} + a(X,L)L$.  In this situation the pushforward map $f_{*}: N_{1}(Y) \to N_{1}(X)$ takes $F_{Y}$ into $F_{X}$.  We say that $f$ is a face contracting morphism if $f_{*}: F_{Y} \to F_{X}$ is not injective.
\end{defi}

One should think of the face contraction condition as a refinement of the $b$-invariant.  Indeed, since the dimensions of $F_{Y}$ and $F_{X}$ are $b(F,Y,f^{*}L)$ and $b(F,X,L)$ respectively, when the $b$-value of $Y$ is greater than the $b$-value of $X$ the map is automatically face contracting.  However, the converse is not true (see \cite[Example 3.6]{LT17}).

When constructing the exceptional set in Manin's Conjecture, we include $f(Y(F))$ whenever $f$ is a dominant face contracting map.  When we have an equality of $a,b$-values and $f$ is not face contracting, we do not remove the contributions.

The following statement summarizes the discussion so far.

\begin{conj}
Let $X$ be a smooth projective geometrically integral Fano variety over a number field.  Set $L = -K_{X}$.  Let $Z$ be the union of $f(Y(F))$ where $f: Y \to X$ varies over all thin maps from a smooth projective geometrically integral variety $Y$ such that $(a(X, L), b(F, X, L)) \leq (a(Y, f^*L), b(F, Y, f^*L))$ and either
\begin{enumerate}
\item $\dim(Y) < \dim(X)$, or
\item $\dim(Y) = \dim(X)$ and $\kappa(K_{Y} + a(Y,f^{*}L)f^{*}L) > 0$, or
\item $\dim(Y) = \dim(X)$, $\kappa(K_{Y} + a(Y,f^{*}L)f^{*}L) = 0$ and $f$ is face contracting.
\end{enumerate}
Then $Z$ coincides with the exceptional set in Manin's Conjecture for $X$ with an anticanonical polarization.
\end{conj}

\subsection{General case}

Let $X$ be a smooth projective geometrically rationally connected and geometrically integral variety over a number field $F$ and let $L$ be a big and nef divisor on $X$.  In brief, one can see the general construction as a ``relative version'' of the simpler construction given above for Fano varieties.  It is convenient to make one more definition.

\begin{defi}
Let $X$ be a smooth projective variety defined over a field $F$ of characteristic $0$. Let $L$ be a big and nef $\mathbb Q$-divisor on $X$. We define
\[
d(X, L) = \dim(X) - \kappa(K_{X} + a(X,L)L).
\]
When $X$ is singular, we define $d(X,L)$ by pulling $L$ back to a smooth model.  Note that $d(X,L)$ is invariant under extension of the ground field.
\end{defi}

We now describe the conjectural exceptional set.  Without loss of generality we may assume that the canonical map $\pi : X \dashrightarrow W$ for $K_{X} + a(X, L)L$ is a morphism.  Let $Z_0$ be the set of rational points contained in the union of $\mathbf B_+(L)$ and a proper closed subset $\pi^{-1}S$ where $S \subset W$ is a proper closed subset such that $\pi$ is smooth over the complement of $S$.

Suppose we have a thin map $f: Y \to X$ such that $d(X,L) = d(Y,f^{*}L)$, $a(X,L) = a(Y,f^{*}L)$ and the image of $Y$ is not contained in $\pi^{-1}S$.  In this situation, we say that $f$ is face contracting if it satisfies the condition of Definition \ref{defi: facecontraction}.  Note that even when $f$ is not dominant this is a sensible definition since the image of $F_{Y}$ under $f_{*}$ is still $F_{X}$ under these assumptions.

Next as $f: Y \rightarrow X$ varies over all thin maps such that $Y$ is geometrically integral and smooth, $d(Y,f^{*}L) < d(X,L)$ and 
\[
(a(X, L), b(F, X, L)) \leq (a(Y, f^*L), b(F, Y, f^*L)),
\]
we define the set $Z_1 \subset X(F)$ by
\[
Z_1 = \bigcup_f f(Y(F)) \subset X(F).
\]

Next as $f: Y \rightarrow X$ varies over all thin maps such that $Y$ is geometrically integral and smooth, $d(Y,f^{*}L) = d(X,L)$, and either
\[
(a(X, L), b(F, X, L)) < (a(Y, f^*L), b(F, Y, f^*L)),
\]
or the $a$ and $b$ values are equal and $f$ is face contracting, we define the set $Z_2 \subset X(F)$ by
\[
Z_2 = \bigcup_f f(Y(F)) \subset X(F).
\]

Finally, as $f: Y \rightarrow X$ varies over all thin maps such that $Y$ is geometrically integral and smooth, $d(Y,f^{*}L) > d(X,L)$ and 
\[
(a(X, L), b(F, X, L)) \leq (a(Y, f^*L), b(F, Y, f^*L)),
\]
we define the set $Z_3 \subset X(F)$ by
\[
Z_3 = \bigcup_f f(Y(F)) \subset X(F).
\]

The union $Z := \cup_{i=0}^{3} Z_{i}$ can be viewed as the potential geometric obstructions to the validity of Manin's Conjecture.  \cite{LST18} predicts that $Z$ coincides with the exceptional set in Manin's Conjecture.



\section{Computing the thin set} \label{sect: computation}

The main theorem of \cite{LST18} shows that the ``geometric exceptional set'' $Z$ constructed in the previous section is always contained in a thin set.

\begin{theo}[\cite{LST18} Theorem 3.5] \label{theo: precisetheorem}
Let $X$ be a geometrically uniruled geometrically integral smooth projective variety defined over a number field $F$ and let $L$ be a big and nef divisor on $X$.  The subsets $Z_{0}$, $Z_1$, $Z_2$, and $Z_{3}$ defined above are contained in a thin subset of $X(F)$.
\end{theo}

In the first subsection, we give a step-by-step construction of a thin subset containing the $Z_{i}$.  In the later subsections, we describe the theoretical results underlying the construction.  Readers who are primarily interested in seeing this construction in action may skip the later subsections and progress directly to examples.

\subsection{The exceptional set}

Here is a list of steps which identifies a thin set containing the set of rational point contributions $Z_{0},Z_{1},Z_{2},Z_{3}$ as in Theorem \ref{theo: precisetheorem}.  While the steps are unfortunately somewhat complicated, in examples it is not uncommon for many of these complications to disappear.  Sections \ref{sect: hypersurfaces} - \ref{sect: bhb} are devoted to examples of this construction; Section \ref{sect: bhb} gives a nice example of how to carry out the process.

\begin{enumerate}
\setcounter{enumi}{-1}
\item Define $V$ to be the proper closed subset $\mathbf{B}_{+}(L)$.  We will increase $V$ subsequently.  

\item Find all adjoint rigid subvarieties with $a(Y,L) \geq a(X,L)$.   
By Theorem \ref{theo: babbound} such subvarieties are either contained in $\mathbf{B}_{+}(L)$ or have $L$-degree bounded above by some constant $C$.  In low dimensions the constant $C$ is known explicitly, giving an explicit bound on the degree of subvarieties $Y$ which must be considered.  For subvarieties with small degree usually one must compute the $a$-invariant by hand by identifying the possible singularity types.

\setcounter{enumi_saved}{\value{enumi}}
\end{enumerate}

Step (1) identifies the bounded subset of $\Chow(X)$ parametrizing all adjoint rigid subvarieties with $a(Y,L) \geq a(X,L)$.  These yield a finite collection of families of subvarieties $\mathcal{U}_{i} \to W_{i}$.  For each index $i$, we add on contributions to our thin set as follows.

\begin{enumerate}
\setcounter{enumi}{\value{enumi_saved}}

\item If $\mathcal{U}_{i}$ is not geometrically integral, then the rational points on $\mathcal{U}_{i}$ will all map into a proper closed subset of $X$.  Enlarge $V$ by adding this subset.

\item If $\mathcal{U}_{i}$ is geometrically integral and the map $\mathcal{U}_{i} \to X$ is not dominant, we enlarge $V$ by adding the image of $\mathcal{U}_{i}$.

\item If $\mathcal{U}_{i}$ is geometrically integral and the map $\mathcal{U}_{i} \to X$ is dominant and has degree $\geq 2$, we will include the image of the rational points on $\mathcal{U}_{i}$ in our thin set. 
\setcounter{enumi_saved}{\value{enumi}}
\end{enumerate}

The only remaining case is when $\mathcal{U}_{i}$ is geometrically integral and the map $\mathcal{U}_{i} \to X$ is birational.  (This includes the possibility that $X$ itself is adjoint rigid and $\mathcal{U}_{i} = X$.)  In this case the construction is more complicated, and we carry it out in several steps.

\begin{enumerate}
\setcounter{enumi}{\value{enumi_saved}}

\item Suppose that $\mathcal{U}_{i}$ is geometrically integral and the map $\mathcal{U}_{i} \to X$ is birational.  Make the following modifications to the family.  After resolving and shrinking the base, we may suppose that $\mathcal{U}_{i} \to W_{i}$ is a smooth family of adjoint rigid varieties and that the $a$-invariant of the fibers is constant.  Let $Q_{i} \subset \mathcal{U}_{i}$ denote the union of all codimension $1$ subvarieties $D$ of fibers $F$ such that $a(D,L) > a(F,L)$.  After shrinking the base further, we may suppose that $Q_{i}$ is flat over $W_{i}$ and that the restriction of the map $\mathcal{U}_{i} \to W_{i}$ to the complement of $Q_{i}$ is topologically locally trivial over $\overline{F}$.  We then take a Galois cover $W_{i}' \to W_{i}$ such that the main component $\mathcal U_{i}'$ of $\mathcal{U}_{i} \times_{W_{i}} W_{i}'$ admits a section $W_{i}' \to \mathcal{U}_{i}'$ avoiding the analogously defined subset $Q_{i}'$.
After these changes, the map $\mathcal{U}_{i}' \to X$ will be proper \'etale over an open subset, and we add the complement to $V$.  We also add the image of $Q_{i}'$ in $X$ to $V$.

\item Let $\mathcal{V}_{i}'$ denote the complement of $Q_{i}'$ in $\mathcal{U}_{i}'$.  The section constructed in the previous step gives us a way of identifying the \'etale fundamental groups of $\overline{P} \cap \overline{\mathcal{V}_{i}'}$ as we  vary the fiber $P$ of $\mathcal{U}_{i}' \to W_{i}'$.

Fix any smooth adjoint rigid variety $P$ in our family $\mathcal{U}_{i}'$.  After passing to the algebraic closure of the ground field, there are up to birational equivalence only finitely many dominant generically finite maps $f: \overline{Y} \to \overline{P}$ such that $a(\overline{Y},f^{*}L) = a(\overline{P},L)$ and $(\overline{Y},f^{*}L)$ is adjoint rigid.  Associated to each such cover there is a finite index subgroup of $\pi_{1}^{\et}(\overline{\mathcal{V}_{i}'} \cap \overline{P})$.  Identify all such finite index subgroups as we vary the fiber $P$ in our family $\mathcal{U}_{i}'$.

To achieve this classification, first note that the degree of $f$ is bounded above by Theorem \ref{theo: babbound}.  Furthermore, by \cite{Sen17} any component of the branch divisor for such a map $f: \overline{Y} \to \overline{P}$ must have $a$-invariant which is strictly larger than the $a$-invariant for $\overline{P}$.  Applying \cite{HJ16} we see there is a finite list of possible branch divisors for $f$.  Together these results show that for any fixed fiber $\overline{P}$ there is a finite set of possible maps $f: \overline{Y} \to \overline{P}$ as above.  Since the degree of the cover is bounded, we see there are only finitely many possible finite index subgroups as we vary $P$. 

\item Using the section constructed above we have an identification
\begin{equation*}
\pi_{1}^{\et}(\overline{\mathcal{V}_{i}'}) = \pi_{1}^{\et}(\overline{\mathcal{V}_{i}'} \cap \overline{P}) \rtimes \pi_{1}^{\et}(\overline{W_{i}}')
\end{equation*}
where $\overline{P}$ is any fiber of $\overline{\mathcal{U}_{i}}' \to \overline{W_{i}}'$.  In step (6) we constructed a finite set of covers of various $\overline{P}$ which yield a finite set of subgroups $\{ \Xi_{i,j} \}_{j \in J}$ of $\pi_{1}^{\et}(\overline{\mathcal{V}_{i}'} \cap \overline{P})$.  For each subgroup $\Xi_{i,j}$, we let $N_{i,j}$ denote the largest finite index subgroup of $\pi_{1}^{\et}(\overline{W_{i}})$ which normalizes $\Xi_{i,j}$ and such that the monodromy action on the cosets of the intersection of conjugates of $\Xi_{i,j}$ is trivial.  The finite index subgroup $\Xi_{i,j} \rtimes N_{i,j}$ of the fundamental group defines a geometric \'etale cover, and after taking completions we obtain a cover $\overline{\mathcal{Y}_{i,j}} \to \overline{\mathcal{U}_{i}'}$ and a family structure $\overline{\mathcal{Y}_{i,j}} \to \overline{T}_{i,j}$ compatible with the map $\overline{\mathcal{U}_{i}'} \to \overline{W}_{i}'$.  We now divide into three cases:
\begin{enumerate}
\item Suppose that the general fiber of $\overline{\mathcal{Y}_{i,j}} \to \overline{T}_{i,j}$ has smaller $a$-value than the general fiber of $\overline{\mathcal{U}_{i}'} \to \overline{W}_{i}'$ or fails to be adjoint rigid.  Then there is a non-dense subset of $\overline{T}_{i,j}$ which parametrizes all fibers with the same $a$-value and are adjoint rigid.  The closure of the image in $X$ of the union of such fibers is a proper subset, and we increase $V$ by adding this set and its Galois conjugates.
\item Otherwise, the general fiber has the same $a$-value and is adjoint rigid. Now suppose that $\overline{\mathcal{Y}_{i,j}} \to \overline{\mathcal{U}_{i}'}$ descends to our ground field $F$ to yield a map $\mathcal{Y}_{i,j} \to \mathcal{U}_{i}'$ in such a way that $\mathcal{Y}_{i,j}$ admits a rational point whose image in $X$ is not contained in $V$ (as constructed so far).  Then we retain the map $\mathcal{Y}_{i,j} \to T_{i,j}$ for the later constructions.  Up to twists and birational equivalence, there are only a finite set of such maps $\mathcal{Y}_{i,j} \to T_{i,j}$.
\item Finally, it may happen that a general fiber has the same $a$-value and is adjoint rigid, but the map $\overline{\mathcal{Y}_{i,j}} \to \overline{\mathcal{U}_{i}'}$ does not descend in such a way as to admit a rational point as explained above.  Let $\overline{\mathcal{Y}_{i,j}}^{\circ}$ denote the open subset consisting of the fibers of the map to $T_{i,j}$ whose $a$-values agree with the generic value. In this case we increase $V$ by adding all Galois conjugates of the proper closed subset of $X$ which is the union of the image of $\overline{\mathcal{Y}_{i,j}} \backslash \overline{\mathcal{Y}_{i,j}}^{\circ}$ with the locus where $\overline{\mathcal{Y}_{i,j}} \to \overline{X}$ is not proper \'etale.
\end{enumerate}

\item At the end of step (7) we have constructed a finite set of families $\mathcal{Y}_{i,j} \to T_{i,j}$.  We next modify each family.  After a birational change we may suppose that $\mathcal{Y}_{i,j}$ is smooth.  After shrinking $T_{i,j}$ we may suppose that the map $T_{i,j} \to W_{i}'$ is \'etale and that the map is smooth over $T_{i,j}$, the $a$-constant is constant over $T_{i,j}$, and the set $Q_{Y}$ defined as before is flat over $T_{i,j}$ and $\overline{\mathcal{Y}_{i,j}} \backslash Q_{Y}$ is locally topologically trivial over $T_{i,j}$.  After a base change we may suppose that $T_{i,j} \to W_{i}'$ is Galois and that the geometric monodromy action on the geometric N\'eron-Severi space of the fibers is trivial.   After another birational change we may suppose that the automorphism group of $\overline{\mathcal{Y}_{i,j}}/\overline{\mathcal{U}}_{i}$ is equal to the full birational automorphism group.  After these changes, the map $\mathcal{Y}_{i,j} \to X$ will be proper \'etale over an open subset, and we add the complement to $V$.  We then replace $\mathcal{Y}_{i,j}$ and $T_{i,j}$ be projective closures while preserving the existence of a map $\mathcal{Y}_{i,j} \to T_{i,j}$.

\item We have a proper closed subset $V \subsetneq X$ from Steps (0)-(8), a finite collection of dominant maps $f_{i}: \mathcal{U}_{i} \to X$ of degree $\geq 2$ from step (4), and a finite collection of dominant maps $g_{i,j}: \mathcal{Y}_{i,j} \to X$ from step (8).  Our thin set is
\begin{equation*}
V(F) \cup \left( \bigcup_{i} f(\mathcal{U}_{i}(F)) \right) \cup  \bigcup_{i,j} \left( \bigcup_{\sigma \in H^{1}(\Gal(\overline{F}/F),\Aut(\overline{Y}_{i,j}/\overline{X}))} g_{i,j}^{\sigma}(\mathcal{Y}_{i,j}^{\sigma}(F)) \right)
\end{equation*}
where the latter union takes place over all twists $\sigma$ such that
\begin{equation*}
(a(X, L), b(F, X, L)) \leq (a(\mathcal{Y}_{i,j}^{\sigma}, g_{i,j}^{\sigma *}L), b(F, \mathcal{Y}_{i,j}^{\sigma}, g_{i,j}^{\sigma*}L)),
\end{equation*}
and $g_{i,j}^{\sigma}$ is face contracting.
\end{enumerate}



\begin{rema}
\cite{LTDuke} conjectured that the covers as in Step (6) must be birationally equivalent to maps that are \'etale in codimension 1.  If true, then one could replace the indirect classification of Step (6) by a calculation involving the \'etale fundamental group.  See Question \ref{ques: etalecodim1} for a further discussion.
\end{rema}

\begin{rema} \label{rema: precisebound}
It is possible that our constructed thin set $Z$ includes more than we should.  We have made no attempt to minimize the closed set $V$ so it might be larger than necessary.  There is an additional issue in Step (4), where we have made no attempt to compare $a,b$-values for the family we remove.  Note that it does not suffice to check the $a,b$-values for the family $\mathcal{U}_{i}$; even if the $a,b$-values for the family $\mathcal{U}_{i}$ are lower, there still might be a cover with higher $a,b$-values, and we should remove the points from such covers.  In practice it is usually easy to see when our thin set $Z$ is too large and which parts of the construction should be modified to give a tight bound.

When $(X, L)$ is adjoint rigid, the extra points that we remove will not affect the asymptotic growth rate. Indeed, removing further proper closed subsets won't affect the asymptotic because $Z_1$ includes all contributions from subvarieties with $a, b$-invariants higher than or equal to the invariants of $X$. Also any cover $f : Y \rightarrow X$ inducing a family in Step (4) satisfies $d(Y, f^*L) < d(X, L)$, so $Z_1$ contains all contributions from such $Y$ with the invariants higher than or equal to the invariants of $X$.
\end{rema}

\subsection{Classifying subvarieties with $\geq$ $a$-value}

We now turn to the theoretical results underlying the description above.  The main finiteness result for subvarieties is:

\begin{theo}[\cite{HJ16} Theorem 1.1] \label{theo: HJ}
Let $X$ be a smooth projective uniruled variety over an algebraically closed field and let $L$ be a big and nef divisor on $X$.  There is a closed proper subset $V \subset X$ which contains every subvariety $Y$ satisfying $a(Y,L) > a(X,L)$.  
\end{theo}

In fact, the components of $V$ are themselves subvarieties of higher $a$-value than $X$.

\begin{proof}
We explain the proof carefully for rational curves on $X$.  Recall that the $a$-invariant of any rational curve on $X$ is $2/L \cdot C$.  We know that any curve not contained in $\mathbf{B}_{+}(L)$ will have positive intersection against $L$.  Thus the rational curves on $X$ not contained in $\mathbf{B}_{+}(L)$  with $a(C,L|_{C}) > a(X,L)$ all have bounded degree.  Such curves must lie on a bounded family in $X$.   Lemma \ref{lemm: dominantfamilies} completes the proof.

For subvarieties of higher dimension, the proof is exactly the same, except that one must appeal to Proposition~\ref{prop: coveredbyadjrigid} to show that one only need consider adjoint rigid varieties and then to Theorem \ref{theo: babbound} to bound the degree. 
\end{proof}

While adjoint rigid subvarieties $Y$ satisfying $a(Y,L) = a(X,L)$ will still lie in a bounded family, they no longer need to be contained in a closed subset since Lemma \ref{lemm: dominantfamilies} no longer gives any restriction.  Nevertheless the geometry of such families is very special.

\begin{theo}[\cite{LTDuke} The proof of Proposition 4.14]
\label{theo: genericallyfinite}
Let $X$ be a smooth projective uniruled variety over an algebraically closed field and let $L$ be a big and nef divisor on $X$.  Suppose that $\mathcal{Y} \to T$ is a family of adjoint rigid subvarieties of $X$ whose general member satisfies $a(Y,L) = a(X,L)$.  Suppose also that the family map $\mathcal{Y} \to X$ is dominant.  If the induced rational map $T \dashrightarrow \mathrm{Chow}(X)$ is generically finite then the map $\mathcal{Y} \to X$ is also generically finite.
\end{theo}

\subsection{Classifying covers with $\geq$ $a$-value}

We now turn to the other extreme and suppose that $f: Y \to X$ is a dominant generically finite morphism of projective varieties.  Lemma \ref{lemm: riemannhurwitz} shows that the $a$-invariant of a cover can only decrease, and we would like to see when equality is attained.  In general there can be infinitely many such maps.

\begin{exam} \label{exam: basechanging}
Let $X$ be $\mathbb{P}^{1} \times \mathbb{P}^{1}$ and let $L = -K_{X}$.  Then for any map of curves $g: C \to \mathbb{P}^{1}$ the induced map $f: C \times \mathbb{P}^{1} \to X$ has that $a(C \times \mathbb{P}^{1},f^{*}L) = 2 = a(X,L)$.
\end{exam}

However, we can expect a better situation when our varieties are adjoint rigid.  The essential ingredient is the following result of \cite{Sen17}.

\begin{theo}[\cite{Sen17} Corollary 2.8] \label{theo: akash}
Let $X$ be a smooth projective geometrically uniruled variety and let $L$ be a big and nef divisor on $X$.  Suppose that $f: Y \to X$ is a dominant generically finite map of smooth varieties such that $(Y,f^{*}L)$ is adjoint rigid and $a(Y,f^{*}L) = a(X,L)$.  Then every component $B$ of the branch divisor on $X$ satisfies $a(B,L) > a(X,L)$.
\end{theo}

By Theorem \ref{theo: HJ} there are only finitely many possible divisors in the branch locus of such a map, and by Theorem~\ref{theo: babbound} the degrees of such maps are bounded.  This yields:

\begin{coro}[\cite{Sen17} Theorem 1.1] \label{coro: akash}
Let $X$ be a smooth projective geometrically uniruled variety and let $L$ be a big and nef divisor on $X$.  Up to birational equivalence, there are only finitely many dominant generically finite maps of smooth varieties $f: Y \to X$ such that $(Y,f^{*}L)$ is adjoint rigid and $a(Y,f^{*}L) = a(X,L)$.
\end{coro}

\subsection{Classifying thin maps with $\geq$ $a$-values}

\begin{defi}
Let $X$ be a smooth projective variety over a field of characteristic $0$ and let $L$ be a big and nef $\mathbb{Q}$-divisor on $X$.  Let $\pi: X \dashrightarrow Z$ be the Iitaka fibration associated to $K_{X} + a(X,L)L$.  An Iitaka base change of $X$ is the normalization of a projective closure of the main component of $T \times_{Z} X$ for some dominant morphism $T \to Z$.
\end{defi}

Just as in Example \ref{exam: basechanging} we can not expect there to be only finitely many thin maps with larger $a$ and $b$ invariants.  Indeed, any Iitaka base change of such a map will also have larger $a$ and $b$ invariants.  However, the following statement proves finiteness up to birational equivalence and Iitaka base change.

\begin{theo}[\cite{LST18} Theorem 4.10] \label{theo: mainfiniteness}
Let $X$ be a uniruled smooth projective variety over an algebraically closed field and let $L$ be a big and nef divisor on $X$.  There is a finite set of thin maps $\{ f_{i}: Y_{i} \to X\}$ with $a(Y_{i},f_{i}^{*}L) \geq a(X,L)$ satisfying the following property.  For any thin map $f: Y \to X$ satisfying $a(Y,f^{*}L) \geq a(X, L)$, after an Iitaka base change to obtain a variety $\widetilde{Y}$ the induced map $\widetilde{f}: \widetilde{Y} \to X$ will factor rationally through some $f_{j}$.  Furthermore, we have
\begin{equation*}
a(Y,f^{*}L) = a(\widetilde{Y}, \widetilde{f}^*L) \leq a(Y_j, f_j^*L)
\end{equation*}
and if equality of $a$-values is achieved then
\begin{equation*}
b(Y,f^{*}L) \leq b(\widetilde{Y},\widetilde{f}^{*}L) \leq b(Y_{j},f_{j}^{*}L).
\end{equation*}
\end{theo}

\subsection{Twists}

In this section, we assume that our ground field is a number field $F$.
Let $(X, L)$ be an adjoint rigid pair defined over $F$.
Using the boundedness of singular Fano varieties, we obtain a boundedness statement over the algebraic closure $\overline{F}$ for adjoint rigid covers $\overline{f} : \overline{Y} \rightarrow \overline{X}$ with the same $a$-invariant as $X$. However, the corresponding statement over the ground field is not true due to the existence of twists. 

Suppose that we have a generically finite dominant $F$-morphism $f : Y \rightarrow X$ of projective varieties defined over $F$. A twist of $f$ is another $F$-morphism $f' : Y' \rightarrow X$ such that over the algebraic closure the two covers are isomorphic, i.e., there exists an $\overline{X}$-isomorphism $\phi : \overline{Y} \cong \overline{Y}'$. An excellent reference about twists is Poonen's book \cite{Poo17}. One can show that the isomorphism classes of twists of $f$ are parameterized by the Galois cohomology
\[
H^1(\mathrm{Gal}(\overline{F}/F), \mathrm{Aut}(\overline{Y}/\overline{X})).
\]

Let us state this one-to-one correspondence explicitly. Suppose that we have a twist $f' : Y' \rightarrow X$ of $f$. We assume that this becomes isomorphic to $f$ over a finite extension $F'/F$. Thus we have $X_{F'}$-isomorphism $\phi : Y_{F'} \cong Y'_{F'}$. Then the associated $1$-cocycle is given by 
\[
\mathrm{Gal}(F'/F) \ni s \mapsto \phi^{-1} \circ \phi^s \in \mathrm{Aut}(\overline{Y}/\overline{X}),
\]
where $\phi^s$ is the conjugate of $\phi$ by $s$.

On the other hand, if we have a $1$-cycle $\mathrm{Gal}(F'/F)\ni s \mapsto \sigma_s \in \mathrm{Aut}(Y_{F'}/X_{F'})$, then we consider the action of $s \in \mathrm{Gal}(F'/F)$ on $Y\otimes F'$ by the composition of $\sigma_s$ and $\mathrm{id}\otimes s$; then we take the quotient
\[
Y^\sigma = (Y\otimes F')/\mathrm{Gal}(F'/F).
\]
This is called as the twist of $Y$ by $\sigma$.

One of the main ingredients of Theorem~\ref{theo: precisetheorem} is the Hilbert Irreducibility Theorem and the following consequence:
\begin{theo}[Hilbert Irreducibility Theorem]
Suppose that $f: Y \rightarrow X$ is a surjective generically finite morphism of normal geometrically integral projective varieties defined over a number field $F$. Suppose that the extension of the geometric function fields $\overline{F}(\overline{Y})/\overline{F}(\overline{X})$ is Galois with Galois group $G$. Then there is a thin set of points $Z \subset X(F)$ such that if $P \in X(F)\setminus Z$, then $f^{-1}(P)$ is irreducible and the corresponding extension of residue fields is Galois with Galois group $G$.
\end{theo}

The following theorem is a consequence of the Hilbert Irreducibility Theorem:
\begin{theo}{\cite[Theorem 5.2]{LST18}}
Let $X$ be a geometrically uniruled normal projective variety defined over a number field $F$. Let $L$ be a big and nef divisor on $X$. We assume that $f :  Y \rightarrow X$ is a dominant generically finite morphism from a normal projective variety. As $\sigma$ varies over all $\sigma \in H^1(F, \mathrm{Aut}(\overline{Y}/\overline{X}))$ such that $Y^\sigma$ is irreducible, 
\[
(a(X, L), b(F, X, L)) \leq (a(Y^\sigma, (f^\sigma)^*L), b(F, Y^\sigma, (f^\sigma)^*L))
\]
and $f^\sigma$ is face contracting, the set 
\[
Z = \bigcup_\sigma f^\sigma(Y^\sigma(F)) \subset X(F)
\]
is contained in a thin set of $X(F)$.
\end{theo}

By combining this theorem with the geometric finiteness in Theorem \ref{theo: mainfiniteness} one can prove that the proposed exceptional set of \cite{LST18} is contained in a thin set.

\section{Hypersurfaces} \label{sect: hypersurfaces}

In the remainder of the paper we will give a number of examples illustrating the computation of thin sets as above.  The easiest examples are varieties with large $a$-invariants.  In this section we discuss low degree hypersurfaces.  Throughout $L$ will denote the restriction to $X$ of the hyperplane class on the ambient projective space.

\begin{exam} \label{exam: projectivespace}
Consider $\mathbb{P}^{n}$ with the hyperplane class $L$.  We have $a(\mathbb{P}^{n},L) = n+1$.  By Lemma \ref{lemm: siu} it is impossible for any variety $Y$ with $\dim(Y)<n$ to have $a(Y,L) \geq a(\mathbb{P}^{n},L)$.  We claim there is also no dominant thin map $f: Y \to \mathbb{P}^{n}$ with $a(Y,f^{*}L) = a(\mathbb{P}^{n},L)$.  Indeed, suppose we have such a map.  Theorem \ref{theo: akash} shows that $(Y,f^{*}L)$ can not be adjoint rigid.  But if $(Y,f^{*}L)$ is not adjoint rigid then it is covered by subvarieties with the same $a$-value, again a contradiction to Lemma \ref{lemm: siu}.

In sum, we see that there are no thin maps $f: Y \to \mathbb{P}^{n}$ with $a(Y,L) \geq a(\mathbb{P}^{n},L)$.  This means that the set of rational points we should remove is the empty set.  For projective space the asymptotic formula in Manin's Conjecture was first confirmed in \cite{Sch79} and was generalized to arbitrary metrization in \cite{Peyre}.
\end{exam}

\begin{exam} \label{exam: quadrics}
Let $Q$ denote a smooth quadric hypersurface of dimension $n \geq 3$.  We have $a(Q,L) = n$.  By Lemma \ref{lemm: siu} it is impossible for any variety $Y$ with $\dim(Y)<n-1$ to have $a(Y,L) \geq a(Q,L)$.  Furthermore, by Proposition \ref{prop: fujitaclassification} any divisor with $a(Y,L) \geq a(Q,L)$ must have $L$-degree $1$, an impossibility in our dimension range.  Using similar logic as in Example \ref{exam: projectivespace} we can use Theorem \ref{theo: akash} to prove there are no dominant thin maps with $a(Y,L) \geq a(X,L)$.

In sum, we see that there are no thin maps $f: Y \to X$ with $a(Y,L) \geq a(X,L)$.  This means that the set we should remove is the empty set. Manin's Conjecture has been confirmed in \cite{FMT89}.
\end{exam}

\begin{exam}
Let $X$ denote a smooth cubic hypersurface of dimension $n \geq 5$.  We have $a(X,L) = n-1$.  By Lemma \ref{lemm: siu} it is impossible for any variety $Y$ with $\dim(Y)<n-2$ to have $a(Y,L) \geq a(X,L)$.  Using Proposition \ref{prop: fujitaclassification} one can  prove that there are no subvarieties with $a(Y,L) \geq a(X,L)$.  Using similar logic as in Example \ref{exam: projectivespace} we can use Theorem \ref{theo: akash} to prove there are no dominant thin maps with $a(Y,L) \geq a(X,L)$.

In sum, we see that there are no thin maps $f: Y \to X$ with $a(Y,L) \geq a(X,L)$.  This means that the set we should remove is the empty set.
Manin's Conjecture is known for $n \geq 7$ by \cite{Hoo94}. For the diagonal cubic we have a result for $n = 6$ by the work of \cite{Vau86}.
\end{exam}

 An interesting but challenging problem is to extend the result to hypersurfaces of any degree.  Based on the properties of rational curves on low degree hypersurfaces proved by \cite{HRS04}, \cite{BK13}, \cite{RY16}, \cite{BV17} and the connection with $a$ and $b$ invariants proved in \cite{LT17}, we expect there are no thin maps $f: Y \to X$ with $\geq a,b$-invariants when $X$ is a general smooth hypersurface in $\mathbb{P}^{n}$ of degree $\leq n-2$ or when $X$ is an arbitrary smooth hypersurface in $\mathbb{P}^{n}$ of degree $\leq n/2$.

\section{Varieties with large group actions} \label{sect: groupcomp}

\begin{exam}[Homogeneous spaces]
Suppose that $X$ is a homogeneous space $G/P$ and $L$ is a big and nef divisor on $X$.  Using the group action, we see that every subvariety of $X$ can be realized as a general member of a dominant family of subvarieties.  By Lemma \ref{lemm: dominantfamilies} we see that $a(Y,L) \leq a(X,L)$ for every subvariety $Y$.  By Theorem \ref{theo: akash} there are no adjoint rigid covers $(Y,f^{*}L)$ with the same $a$-value as $X$. 
Any fibration on $X$ has the form $f : G/P \rightarrow G/M$ where $M \supset P$ is a parabolic subgroup of $G$. Any adjoint rigid subvariety $Y$ with $a(Y, L)= a(X, L)$ is a fiber of some such fibration $f$.
Indeed, suppose that $Y$ contains a point coming from the identity and let $N$ be the stabilizer of $Y$, i.e., $N = \{g \in G \mid g \cdot Y = Y\}$.
Then translations give us an injection $G/N \rightarrow \mathrm{Chow}(X)$. Let $\mathcal U$ be the pullback of the universal family on $G/N$ which is a family of adjoint rigid subvarieties. Then we have a dominant morphism
\[
\mathcal U \rightarrow X
\]
Since the induced map $G/N \rightarrow \mathrm{Chow}(X)$ is injective, it follows from Theorem~\ref{theo: genericallyfinite} that the above map is generically finite. In particular, we conclude that $\dim Y = \dim N - \dim P$. Since we have the injection $N/N\cap P \hookrightarrow Y$, this implies that $N \supset P$ and $Y = N/P$. Thus our assertion follows.

Suppose that $L$ is proportional to $-K_X$. Then there is no thin map $g : Y \rightarrow X$ such that 
\[
(a(X, L), b(F, X, L)) \leq (a(Y, g^*L), b(F, Y, g^*L)).
\]
Indeed, without loss of generality we may assume that $(Y, g^*L)$ is adjoint rigid. Then $g(Y)$ is an adjoint rigid subvariety so the above discussion shows that $g(Y)$ itself is a homogenous space. This shows that $g(Y)$ admits no adjoint rigid cover. Thus $Y$ is an adjoint rigid subvariety of $X$ and it is smooth. It follows from \cite[Proposition 3.5]{HTT15} that we have {\it the balanced property}, i.e.,
\[
(a(Y, L), b(F, Y, L)) < (a(X, L), b(F, X, L)).
\]
Thus our assertion follows.

In sum we conclude that our proposed exceptional set is the empty set. Another direct approach to this problem is discussed in \cite[Example 3.9]{LST18}.
Manin's Conjecture for generalized flag varieties has been confirmed in \cite{FMT89}.
\end{exam}

\begin{exam}[Equivariant compactifications of vector groups]
Let $G = \mathbb G_a^n$ be the $n$-dimensional vector space and let $X$ be a smooth projective equivariant compactification of $G$ defined over a number field $F$. Manin's Conjecture has been proved after removing rational points on the boundary in \cite{CLT02}. Let $L$ be a big and nef divisor on $X$. It is clear that any subvariety with higher $a$-invariant is contained in the boundary. Since any vector group is simply connected, we conclude that there is no adjoint rigid cover with the same $a$-invariant as $X$. Let $Y$ be an adjoint rigid subvariety with the same $a$-invariant as $X$ and  meeting with $G$. Let $N$ be the stabilizer of $Y$. We fix a section $G/N \rightarrow G$. Then using this section translations define a dominant morphism
\[
G/N \times Y \rightarrow X.
\]
By Theorem~\ref{theo: genericallyfinite}, this map is generically finite. In particular, we conclude that $\dim Y = \dim N$. Since we have the injection $N \hookrightarrow Y$, this implies that $Y$ is an equivariant compactification of $N$. When $(X, L)$ is adjoint rigid, it follows from \cite[Theorem 5.4]{HTT15} that we have
\[
(a(Y, L), b(F, Y, L)) < (a(X, L), b(F, X, L)).
\]
Hence we conclude that our proposed exceptional set is contained in the boundary. Readers can find a more arithmetic treatment of this problem in \cite[Example 3.10]{LST18}.
\end{exam}

\begin{exam}[Toric varieties] \label{exam: toricvarieties}
Let $X$ be a geometrically integral smooth toric variety defined over a number field $F$ and let $L$ be a big and nef divisor on $X$. Manin's Conjecture for such a variety was proved in \cite{BT-general}, \cite{BT-0}, and \cite{Sal98} after removing rational points on the boundary.

The maps $f: Y \to X$ such that $a(Y,f^{*}L) \geq a(X,L)$ are controlled by the following restrictions:
\begin{itemize}
\item By Lemma \ref{lemm: dominantfamilies} any subvariety with larger $a$-value must be contained in the toric boundary.
\item If we have a dominant family of subvarieties $Y$ with $a(Y,L) = a(X,L)$, then by the argument of \cite[Theorem 9.5]{LTDuke} the varieties $Y$ must be the fibers of a rational toric map $X \dashrightarrow W$.
\item Let $f: Y \to X$ be a dominant thin map such that $a(Y,f^{*}L) = a(X,L)$ and $(Y,L)$ is adjoint rigid.  We claim that after applying some birational modification, $Y$ is toric and $f$ is equivariant. Indeed, taking the Stein factorization, we may assume that $f$ is finite. By Theorem~\ref{theo: akash}, the branch locus is contained in the boundary, so $Y$ contains a torus $T$ as a dense open set. For any $t \in T$, we may define a birational map $g_t : Y \dashrightarrow Y$ by translation on $T$. We also have the corresponding automorphism $h_t : X \rightarrow X$. Let $A$ be an ample divisor on $X$. Then we see that $g_t^*f^*A \sim f^*A$. Thus $g_t$ must be regular.
By applying an equivariant resolution we may assume that $Y$ is a smooth toric variety and $f$ is equivariant.  

We next claim that $f$ is birational to a map $f': Y' \dashrightarrow X'$ that is \'etale in codimension $1$.  It suffices to prove that the branch divisor $B$ on $X$ has Iitaka dimension $0$.  One can then contract $B$ by a toric map $X \dashrightarrow X'$, and the Stein factorization of the induced map $Y \dashrightarrow X'$ is the map $f'$.  To prove the claim on the branch divisor $B$, we note that there is some multiple $mR$ of the ramification divisor $R$ on $Y$ such that $mR - f^{*}B$ is pseudo-effective.  Indeed, it is clear that every $f$-exceptional divisor has positive coefficient in $R$.  Furthermore, the map $f_{*}$ induces a bijection between the non-$f$-exceptional torus invariant divisors on $Y$ and the torus-invariant divisors on $X$.  In particular, the pullback of any component of the branch divisor is a positive linear combination of $f$-exceptional divisors and the unique torus-invariant divisor lying above it.  Since $mR-f^{*}B$ is pseudo-effective and $R$ has Iitaka dimension $0$, $B$ also must have Iitaka dimension $0$.
\end{itemize}
Together, these results show that the only thin maps $f: Y \to X$ with $a(Y,f^{*}L) \geq a(X,L)$ must come from (compositions of) maps of toric varieties and inclusions of torus-invariant subvarieties.  In principle one should be able to classify all such maps combinatorially.

In order to show that the proposed exceptional set of \cite{LST18} is contained in the boundary, it suffices to prove that there is no dominant thin map $f: Y \to X$ such that
\begin{equation*}
(a(X,L),b(F,X,L)) \leq (a(Y,f^{*}L),b(F,Y,f^{*}L))
\end{equation*}
and $f$ is face contracting.  Suppose $f$ is a dominant thin map such that $a(Y,f^{*}L) = a(X,L)$.  By the reasoning above we may assume that $f$ is toric.  Write
\begin{equation*}
K_{Y} + a(Y,f^{*}L)f^{*}L = f^{*}(K_{X} + a(X,L)L) + R.
\end{equation*}
Then $R$ must include every $f$-exceptional divisor in its support.  Furthermore, the map $f_{*}$ induces a bijection between the non-$f$-exceptional torus invariant divisors on $Y$ and the torus-invariant divisors on $X$ so that there is a unique torus invariant divisor on $Y$ which dominates a given toric invariant divisor on $X$.  Thus we see that the map $f$ can never be face contracting.  This shows that our proposed exceptional set is compatible with the known results on Manin's Conjecture.  An alternative arithmetic approach to this problem is given in \cite[Example 3.10]{LST18}.
\end{exam}

One can also consider equivariant compactifications of other linear algebraic groups. Manin's Conjecture for such varieties has been proved in \cite{STBT07} and \cite{ST16}.  One can show that our proposed exceptional set is contained in the boundary using the arithmetic arguments of \cite[Example 3.10]{LST18}.

\section{Surfaces} \label{sect: surfaces}

Manin's Conjecture has been verified for many surfaces; see \cite{dBBD07}, \cite{Bro09}, \cite{Bro10}, and \cite{dBBP12}.  In all known examples the exceptional set is not Zariski dense.  This coheres with the following result:

\begin{theo} \label{theo: delpezzosurfaces}
Let $X$ be a geometrically integral smooth del Pezzo surface over a number field and let $L = -K_{X}$.  If $X$ has degree $1$ then suppose furthermore that either $\rho(X) \geq 2$ or that $\overline{X}$ is general in moduli.  Then $\cup_{i=0}^{3} Z_{i}$ is contained in a proper closed subset of $X$.  More precisely, it consists of the rational points on the finite set of irreducible curves $C$ that are defined over the ground field and satisfy either
\begin{itemize}
\item $L \cdot C = 1$, or
\item $L \cdot C = 2$ and $C$ is singular.
\end{itemize}
\end{theo}

In this section we will let $X$ denote a smooth geometrically integral geometrically uniruled projective surface defined over a number field and let $L$ denote a big and nef divisor on $X$.  We divide into two cases based on the Iitaka dimension of $K_{X} + a(X,
L)L$.

\subsection{Iitaka dimension $1$}

The easier case is when $\kappa(K_{X} + a(X,L)L) = 1$.  In this case $X$ admits a morphism $\pi: X \to Z$ whose fibers are rational curves with the same $a$-values as $X$.  The only subvarieties with $a(Y,L) \geq a(X,L)$ will be components of the fibers of $\pi$ or rational curves satisfying $L \cdot C = 0$.  To construct the thin set $Z$ we will just need to remove the rational points on $\mathbf{B}_{+}(L)$ and in the locus where $\pi$ is not smooth.

\subsection{Iitaka dimension $0$}

Suppose that $(X,L)$ is adjoint rigid.  The first subsection below discusses the behavior of the $a$-invariant for subvarieties and covers.  The second subsection carries out Steps (0)-(9) of Section \ref{sect: computation} using the techniques of \cite{LTDuke}.  Since the behavior of $a$ and $b$ invariants for surfaces over an algebraically closed field is completely understood, the challenge is to understand what happens over a number field.

\subsubsection{Classifying thin maps}  
When $(X,L)$ is adjoint rigid, we can run the minimal model program for $K_{X} + a(X,L)L$ by sequentially contracting a sequence of Galois orbits of geometric $(-1)$-curves which are disjoint.  This yields a morphism $\phi: X \to X'$ to a weak del Pezzo surface $X'$ whose Picard rank is the $b$-value of $X$.  Most of the important features of the $a$ and $b$ invariants on $X$ are best computed by referring to this surface $X'$.

We start by discussing the set of thin maps $f: Y \to X$ where the $a$-invariants for $Y$ are at least as large as those of $X$.  The only possible adjoint rigid subvarieties of $X$ satisfying this condition on $a$-invariants are rational curves $C$ with $L \cdot C \leq 2/a(X,L)$.  There will be a finite set of such curves $C$ which do not deform.   There will also be a finite set of dominant families, which are characterized by:

\begin{lemm} \label{lemm: curvesarefibers}
Any dominant family of rational curves whose general member $C$ satisfies $a(C,L) = a(X,L)$ will lie in one of the following categories.
\begin{enumerate}
\item The curves $C$ are the fibers of a morphism $f: X \to T$.
\item The minimal model of $K_{X} + a(X,L)L$ is a weak del Pezzo surface $X'$ of degree $2$ and $C$ is the strict transform of a rational curve in $|-K_{X'}|$.
\item The minimal model of $K_{X } +a(X,L)L$ is a weak del Pezzo surface $X'$ of degree $1$ and $C$ is the strict transform of a rational curve in $|-2K_{X'}|$.
\item The minimal model of $K_{X } +a(X,L)L$ is a weak del Pezzo surface $X'$ of degree $1$ and $C$ is the strict transform of a rational curve $C'$ on $X'$ satisfying $(C')^{2} = 2$, $-K_{X'}.C' = 2$.
\end{enumerate}
Note that any rational curve in (4) is a member of the linear system of the pullback of the anticanonical divisor of a weak del Pezzo surface of degree $2$.
\end{lemm}

\begin{proof}
By \cite[Proposition 4.14]{LTDuke} such curves must form a $1$-dimensional family.  Using the deformation theory of rational curves we see that $-K_{X} \cdot C = 2$.  Since $a(X,L)L \cdot C = 2$, we see that $(K_{X} + a(X,L)L) \cdot C = 0$.  Let $\phi: X \to X'$ be the result of the $(K_{X} + a(X,L)L)$-MMP.  Then we see that the class of $C$ is pulled back from the weak del Pezzo surface $X'$.  Letting $C'$ denote the image of the curves $C$ on $X'$, we still have $-K_{X'} \cdot C' = 2$.

Since $-K_{X'} \cdot C' = 2$, by adjunction $C'^2$ must be even and non-negative. By the Hodge index theorem we must have $(-K_{X'})^2(C')^2 - 4 \leq 0$.  This gives four possibilities for $K_{X}^{2}$ and $C'^{2}$.  Three of them match the conditions (2)-(4) in the statement above, and all we need to show is that if $C'^{2}=0$ then $C'$ (and hence $C$) is a fiber of a morphism.  But this follows from the Basepoint Free Theorem for the surface $X'$.
\end{proof}

We next discuss adjoint rigid covers of $X$.  Let $X'$ denote the weak del Pezzo surface obtained by running the MMP on $K_{X} + a(X,L)L$.  By further passing to the anticanonical model of $X'$, we find a morphism $\phi: X \to \widetilde{X}$ where $\widetilde{X}$ is a Fano surface with canonical singularities.  Such a surface $\widetilde{X}$ is known as a Gorenstein log del Pezzo surface.

To classify dominant maps with the same $a$-value, it suffices to first work over an algebraically closed field and then to see which ones descend to the ground field.  The first step is accomplished by:

\begin{theo}[\cite{LTDuke} Corollary 6.9] \label{theo: acoversurface} 
Suppose that $F$ is algebraically closed.  Let $\widetilde{X}$ be the Gorenstein log del Pezzo surface constructed above and fix a smooth point $x$.  There is a bijection between finite index subgroups of $\pi_{1}^{\et}(\widetilde{X}^{sm},x)$ and birational equivalence classes of thin maps $f: Y \to X$ such that $a(Y,f^{*}L) = a(X,L)$ and $(Y,f^{*}L)$ is adjoint rigid.
\end{theo}


The bijection associates to a finite index subgroup a projective closure of the corresponding \'etale cover of $\widetilde{X}^{sm}$.  We note for emphasis that the divisor $L$ is \emph{not} necessarily pulled back from $\widetilde{X}$, so that the relationship between the $a$-value of $L$ and covers of $\widetilde{X}$ is a bit surprising.

\begin{rema}
In fact we can say a little bit more.  Suppose again that $F$ is algebraically closed.  In this situation \cite[Lemma 2]{MZ93} shows there is a morphism $g: \widetilde{X} \to W$ where $W$ is a relatively minimal Gorenstein log del Pezzo surface and we have $\pi_{1}^{\et}(\widetilde{X}^{sm},x) = \pi_{1}^{\et}(W^{sm},w)$. The singularity types and quasi-universal covers of the relatively minimal Gorenstein log del Pezzo surfaces $W$ have been classified in \cite{MZ88}, \cite{MZ93}, \cite{ye02}.
\end{rema}

\subsubsection{Constructing a thin set}
We next carry out Steps (0)-(9) of Section \ref{sect: computation} for an adjoint rigid surface pair $(X,L)$.  In step (0), we set $V$ to be the set of rational curves satisfying $L \cdot C = 0$.  Step (1) was carried out over $\overline{F}$ above, and we just need to take those curves which are defined over $F$.  It is helpful to note the following lemma:

\begin{lemm} \label{lemm: surfaceclassification}
If $X$ is adjoint rigid and admits a covering family of curves satisfying $a(C,L) = a(X,L)$ then one of the following conditions holds:
\begin{enumerate}
\item $b(X,L) > 1$.
\item The minimal model of $K_{X} + a(X,L)L$ is a del Pezzo surface $X'$ of Picard rank $1$ and degree $2$ and $C$ is the strict transform of a rational curve in $|-K_{X'}|$.  Such curves are parametrized by an irreducible curve of geometric genus $3$.
\item The minimal model of $K_{X } +a(X,L)L$ is a del Pezzo surface $X'$ of Picard rank $1$ and degree $1$ and $C$ is the strict transform of a rational curve in $|-2K_{X'}|$.  When $\overline{X'}$ is general in moduli such curves are parametrized by an irreducible curve of geometric genus $\geq 2$.
\end{enumerate}
\end{lemm}

\begin{proof}
The classification into three types follows immediately from the argument of Lemma \ref{lemm: curvesarefibers} and the fact that any weak del Pezzo surface of Picard rank $1$ is an honest del Pezzo surface.  In case (2) $X'$ is the double cover of $\mathbb{P}^{2}$ ramified over a smooth quartic, and the rational curves in $|-K_{X'}|$ are parametrized by the curve dual to the quartic which has genus $3$.  In case (3) the linear system $|-2K_{X'}|$ realizes $X'$ as the double cover of a quadric cone.  The ramification locus is the union of the vertex and a smooth genus $4$ curve $T$ cut out by a cubic.  The rational curves in $|-2K_{X'}|$ correspond to the hyperplanes in $\mathbb{P}^{3}$ that are bitangent to $T$.  Let $M$ denote the curve parametrizing these rational curves.  Since $\overline{X'}$ is general \cite[Theorem 21 Claim 11]{Lubbes14} shows that $M$ is irreducible.  Note that $M$ is naturally a subvariety of $\Sym^{2}(T)$.  Since $T$ is canonically embedded, it is not hyperelliptic, thus $\Sym^{2}(T)$ embeds into the Jacobian $\Jac(T)$.  In particular $M$ is not rational, and it can only be elliptic if it is smooth.  But an intersection calculation on $T \times T$ (\cite{Lubbes14}) shows that $M$ has high arithmetic genus, so it can not be an elliptic curve.
\end{proof}

In Step (2) any family of curves which is not geometrically irreducible will only have finitely many members defined over the ground field, and we remove the points on this finite set of curves. (Note that if such a family of curves comes from a fibration $X \to T$ then no smooth fiber can carry any points.  In particular, if $X$ is a del Pezzo surface then the curves whose points we remove must satisfy one of the conditions of Theorem \ref{theo: delpezzosurfaces}.)  In Step (3) we enlarge $V$ to include any curves with $a(C,L) \geq a(X,L)$ which do not deform.  Step (4) applies only in cases (2) and (3) of Lemma \ref{lemm: surfaceclassification}.  However, in each case the family of curves is parametrized by a curve of genus $>1$ and so these rational point contributions are not Zariski dense.

In Step (5), we start with a finite set of fibrations $p_{i}: X \to T_{i}$ whose fibers are rational curves.  We shrink each to the smooth locus and add the complement to $V$.  Since a Hirzebruch surface admits a section, no further changes are necessary.

For the remaining steps, there are two kinds of families to consider.  First, there are the families of rational curves constructed in Step (5).  In Step (6) we must compute all non-trivial $a$-covers of the rational curves, but as we have already seen there are no such covers.  In Steps (7) and (8) we do not make any changes to these families.  The other type of map to consider is the finite set of adjoint rigid dominant covers $f: Y \to X$.  These are classified in some sense by Theorem \ref{theo: acoversurface} accomplishing Step (6).  In Step (7) we must find which ones descend to the ground field and admit a rational point, and in Step (8) we do not make any change to the cover.

Finally, in Step (9) we must calculate which of these families to discount.  First, consider the families of rational curves which are fibers of a morphism.  Here the $b$-value of the curves is $1$ but by Lemma \ref{lemm: surfaceclassification} the $b$-value of $X$ is greater than $1$.  Thus we do not need to remove any point contributions.  Second, consider the adjoint rigid covers of $X$.  Over an algebraically closed field the $b$-value of a cover will always be smaller than the $b$-value of $X$.  However, we do not know whether this is true over a number field.  Thus, we can only complete this step in the situation when $X$ is a smooth del Pezzo surface and $L=-K_{X}$ by using the following result.

\begin{prop}[\cite{LTDuke} Theorem 6.2]
Let $X$ be a smooth del Pezzo surface and $L = -K_{X}$.  Then there is no dominant generically finite map $f: Y \to X$ from a smooth variety $Y$ satisfying $a(Y,f^{*}L) = a(X,L)$ and $(Y,f^{*}L)$ is adjoint rigid.
\end{prop}

Combining all these steps, we see that in the special case when $X$ is a del Pezzo surface and $L = -K_{X}$, the proposed exceptional set of \cite{LST18} is exactly described by Theorem \ref{theo: delpezzosurfaces}.   It seems worthwhile to finish the classification for surfaces by answering the following questions:

\begin{ques} \label{ques: ratcomponents}
Let $X$ be a geometrically integral smooth del Pezzo surface of degree $1$ over a number field.  Let $M$ denote the curve parametrizing rational curves in $|-2K_{X}|$.  Does every component of $M$ have genus $\geq 2$?
\end{ques}

\begin{ques} \label{ques: weakdelpezzo}
Let $X$ be a geometrically integral smooth weak del Pezzo surface over a number field and let $L = -K_{X}$.  Suppose that $f: Y \to X$ is a dominant generically finite map such that $a(Y,f^{*}L) = a(X,L)$ and $(Y,f^{*}L)$ is adjoint rigid.  Then is $b(F,Y,f^{*}L) < b(F,X,L)$?
\end{ques}

If these questions are answered affirmatively, then we expect the exceptional set for a surface to always be contained in a Zariski closed proper subset.

\section{Fano threefolds of Picard rank $1$} \label{sect: fanothreefolds}

We next discuss the $a$ and $b$ invariants for Fano threefolds of Picard rank $1$.  Here we will always let $L$ be a generator of the Picard group so that $a(X,L)$ is the same as the index of $X$ and $b(F, X,L)=1$.

\subsection{Index $2$ Fano threefolds}

\begin{theo}[\cite{LT17} Section 7.2] \label{theo: index2fanothreefolds}
Let $X$ be a Fano threefold of Picard rank $1$ and index $2$ satisfying $L^{3} \geq 2$.  
Then:
\begin{enumerate}
\item There is no subvariety $Y \subset X$ satisfying $a(Y,L) > a(X,L)$.
\item There is no dominant generically finite map $f: Y \to X$ such that $a(Y,L) = a(X,L)$ and $(Y,f^{*}L)$ is adjoint rigid.
\item The only family of adjoint rigid subvarieties with $a(Y,L) = a(X,L)$ is the family of lines.
\end{enumerate}
\end{theo}

Together, these imply the following.  Suppose that $f: Y \to X$ is any thin map satisfying $a(Y,L) = a(X,L)$.  Then $f$ must factor through the universal family of lines.
Let $\pi : \mathcal U \rightarrow W$ be the universal family of lines. Then the proposed exceptional set of \cite{LST18} coincides with $\mathcal U(F)\cup B(F)$ where $B$ is the branch locus of $s : \mathcal U \rightarrow X$.

Over an algebraically closed field, the variety of lines $W$ is of general type for a cubic threefold, an abelian surface for a complete intersection of two quadrics, and $\mathbb P^2$ for $L^3 = 5$. Thus depending on $X$ the Zariski density of rational points on $\mathcal U$ varies. Note that for a cubic threefold, the variety $W$ admits an embedding into the intermediate Jacobian of $X$, in particular it follows from Faltings' theorem that the set of rational points on $W$ is not Zariski dense. On the other hand, in the cases of $L^3 = 4, 5$, $W$ satisfies the potential density.

\begin{exam}
Let $X$ be a smooth complete intersection of two quadric hypersurfaces in $\mathbb{P}^{5}$ and let $L$ denote the hyperplane class on $X$.  Since $a(X,L)=2$ it is easy to see that the only curves with the same $a$-value must be lines, and no curve can have larger $a$-value.

Suppose now that $S$ is a surface in $X$ satisfying $a(S,L) \geq a(X,L)$.  Since $(2L)^{2} \cdot S \geq 16$, Theorem \ref{theo: surfacebabbound} shows that $S$ can not be adjoint rigid.  This means that $S$ must be covered by curves with the same $a$-value as $S$.  In particular, we see that $S$ must be fibered by lines and that $a(S,L) = a(X,L)$.

Since $X$ does not admit any subvarieties with higher $a$-value, Theorem \ref{theo: akash} shows that there is no dominant generically finite map $f: Y \to X$ such that $a(Y,L) = a(X,L)$ and $(Y,f^{*}L)$ is adjoint rigid.
\end{exam}

\subsection{Index $1$ Fano threefolds}

\begin{theo}[\cite{LTDuke} Theorem 1.9 and \cite{LT18} Theorem 5.1]
Let $X$ be a general Fano threefold of Picard rank $1$ and index $1$ such that $-K_{X}$ is very ample.  Then:
\begin{enumerate}
\item Let $D$ denote the divisor in $X$ swept out by the family of lines.  Then $D$ contains every subvariety $Y$ satisfying $a(Y,L) > a(X,L)$.
\item There is no dominant generically finite map $f: Y \to X$ such that $a(Y,L) = a(X,L)$ and $(Y,f^{*}L)$ is adjoint rigid.
\item The only family of adjoint rigid subvarieties with $a(Y,L) = a(X,L)$ is the family of conics.
\end{enumerate}
\end{theo}

This result implies that if $f: Y \to X$ is any thin map satisfying $a(Y,L) \geq a(X,L)$ then $f$ must factor through $D$ or through the universal family of conics.  
Let $\pi : \mathcal U \rightarrow W$ be the universal family of conics. Then the proposed exceptional set of \cite{LST18} is
\[
\mathcal U(F) \cup D(F) \cup B(F)
\]
where $B$ is the branch locus of $s: \mathcal U \rightarrow X$.  When $H^3 = 22, 18$, the variety $W$ satisfies the potential density. When $H^3 = 16, 14, 12$, the set of rational points on $W$ is not Zariski dense. See \cite{KPS18} for more details.

\begin{exam}
Let $X$ denote a general quartic hypersurface in $\mathbb{P}^{4}$ and let $L$ denote the hyperplane class.  The only curves with $a(C,L) \geq 1$ will be the lines and conics on $X$.  Let $D$ denote the surface swept out by lines.

Suppose now that $S$ is a surface in $X$ satisfying $a(S,L) \geq a(X,L)$.  Our goal is to show that $S$ cannot be adjoint rigid.  We know that $S$ will be linearly equivalent to some multiple of $L$.  If $S \sim kL$ for $k \geq 3$ then $L^{2} \cdot S > 9$ and Theorem \ref{theo: surfacebabbound} shows that $S$ can not be adjoint rigid.  The cases where $S \sim L$ or $S \sim 2L$ must be addressed by hand.  \cite{LT18} shows that when $X$ is general any element of $|L|$ is normal and that any element of $|2L|$ is normal after perhaps blowing up a line contained in the surface.  One can then appeal to Lemma \ref{lemm: normalsurfacebound}.

Theorem \ref{theo: akash} states that any dominant generically finite map $f: Y \to X$ such that $a(Y,L) = a(X,L)$ and $(Y,f^{*}L)$ is adjoint rigid must be ramified over the divisor $D$.  It is not immediately clear that no such map exists.  However, \cite[Theorem 1.9]{LTDuke} proves there can be no such cover using a different approach based on Question \ref{ques: etalecodim1}.
\end{exam}

\section{Del Pezzo varieties} \label{sect: delpezzos}

A smooth Del Pezzo variety is a projective Fano variety $X$ of dimension $n$ admitting an ample divisor $L$ such that $K_{X} = -(n-1)L$.  The Del Pezzo varieties are classified in \cite{fujita80}, \cite{fujita81}.  In dimension $\geq 4$ the possible degrees of a Del Pezzo variety are $1 \leq L^{\dim X} \leq 5$, and we will only consider the case when the degree is $\geq 2$.

\begin{theo} \label{theo: delpezzos}
Let $X$ be a general smooth del Pezzo variety of degree $2,3$ and dimension $\geq 4$, of degree $4$ and dimension $\geq 5$, or of degree $5$ and dimension $6$. Then there are no thin maps $f: Y \to X$ such that $a(Y,f^{*}L) \geq a(X,L)$.
\end{theo}

\begin{proof}
Suppose first of all that $X$ is a general smooth del Pezzo variety of degree $2 \leq d \leq 5$ and dimension $\geq 4$.  We would like to classify all adjoint rigid subvarieties $Y \subset X$ satisfying $a(Y,L) \geq a(X,L)$.  It is clear that there are no such curves $Y$ since $a(X,L) \geq 3$.   Suppose we cut down $X$ by a general hyperplane $H$.  Applying Theorem \ref{theo: hyperplanesection}, we see that $a(Y \cap H,L) \geq a(H,L)$.  Applying this construction inductively, we will eventually obtain a threefold $X'$ which is a general smooth Fano threefold of Picard rank $1$, index $2$, and degree $\geq 2$.  As discussed in Theorem \ref{theo: index2fanothreefolds}, the only subvarieties of $X'$ with $\geq$ $a$-value are the covering family of lines.  Thus we see that any subvariety $Y \subset X$ with $\geq$ $a$-value must have codimension $2$ and must be linear in the ambient projective space.

We now give a case-by-case analysis showing that a general $X$ in our dimension range does not contain a codimension $2$ linear space.  In fact, using the same argument cutting down by hypersurfaces, it suffices to prove this when $X$ has the minimal possible dimension for its degree.
\begin{itemize}
\item[] \textbf{Degree 2:} $X$ is a double cover of $\mathbb{P}^{4}$ ramified along a quartic $B$.  The image of any plane contained in $X$ will be a plane $P$ in $\mathbb{P}^{4}$ such that $P \cap B$ is everywhere non-reduced.  A general quartic will not contain any planes.  Consider the incidence correspondence consisting of pairs $(Q,P)$ of a quartic hypersurface $Q$ whose restriction to a plane $P$ is a double conic.  This has dimension $66$, showing that the general quartic in $\mathbb{P}^{4}$ does not restrict to a double conic along any plane.
\item[] \textbf{Degree 3:} $X$ is a cubic fourfold, and a general cubic fourfold does not contain a plane.
\item[] \textbf{Degree 4:} $X$ is a complete intersection of two quadrics in $\mathbb{P}^{7}$.  Such a variety can never contain a $\mathbb{P}^{3}$ by \cite[Corollary 2.4]{Reid}.
\item[] \textbf{Degree 5:} $X$ is $G(2,5)$, so $X$ does not contain any linear spaces of dimension $4$.
\end{itemize}
\end{proof}

We note in passing that the dimension ranges in Theorem \ref{theo: delpezzos} are optimal.  If $X$ is a del Pezzo variety of dimension $4$ and degree $4$, then $X$ is a complete intersection of two quadrics in $\mathbb{P}^{6}$ and such varieties always contain planes by \cite[Corollary 3.3]{Reid}.  If $X$ is a del Pezzo variety of dimension $5$ and degree $5$, then $X$ is a hyperplane section of $G(2,5)$ and such varieties always contain a $\mathbb{P}^{3}$ by \cite[Theorem 1.2]{AC12}.

For many of these varieties Manin's Conjecture is known (with empty exceptional set).  As we discussed earlier, Manin's Conjecture is known for cubic hypersurfaces of dimension $\geq 7$ by \cite{Hoo94} and there are partial results in dimension $6$ by \cite{Vau86}.
Birch's circle method (\cite{Bir61}) can handle a complete intersection of two quadrics when the dimension is $> 9 $. 
Manin's Conjecture is also verified for $G(2,5)$ by the general results of \cite{FMT89}.

\section{Hypersurfaces in biprojective space} \label{sect: bhb}

We work over the ground field $\mathbb{Q}$.  Let $X$ be the smooth hypersurface of degree $(1, 2)$ in the product of projective spaces $\mathbb P^3_x \times \mathbb P^3_y$ defined by the equation $\sum x_{i}y_{i}^{2} = 0$. In this section we compute the exceptional set defined in \cite{LST18} for this example. Manin's Conjecture has been confirmed in \cite{BHB18} after removing a thin exceptional set and we prove that the exceptional set proposed by \cite{LST18} coincides with this exceptional set.

We let $\pi_{1}$ and $\pi_{2}$ denote the two projection maps from $X$ to $\mathbb{P}^{3}$.  Let $h_1$ be the pullback of the hyperplane class via the first projection $\pi_1$ and $h_2$ be the pullback of the hyperplane class via the second projection $\pi_2$.  The lattice $N^{1}(X)_{\mathbb{Z}}$ is generated by $h_{1}$ and $h_{2}$.

We set $L = -K_{X} = 3h_{1} + 2h_{2}$.  Note that $a(X,L) = 1$ and $b(\mathbb{Q},X,L) = 2$.  We start by working over $\overline{\mathbb{Q}}$ and computing the $a$-invariants and $b$-invariants of subvarieties with respect to $L$.  The first step is to understand the singularities of some hyperplane sections of $X$.

\begin{defi}
In this section we will let $S$ denote the union of the singular fibers of $\pi_{1}$. Note that $S$ is the union of the $\pi_{1}$-preimages of the coordinate hyperplanes in $\mathbb{P}^{3}$.
\end{defi}

\begin{lemm} \label{lemm: threefoldquadricbundle}
Let $\ell \subset \mathbb{P}^{3}$ be a line and let $Y$ denote the $\pi_{1}$-preimage of $\ell$ in $X$.  Suppose that $Y \not \subset S$.  Then $Y$ is normal with canonical singularities.
\end{lemm}

\begin{proof}
We first prove that $Y$ is normal.  It suffices to show that the singular locus of $Y$ has codimension $\geq 2$.  Suppose that $y$ is a singular point of $Y$.  After an appropriate coordinate change, we may suppose that in a local $\mathbb{A}^{1}_{s} \times \mathbb{A}^{3}_{y_{1},y_{2},y_{3}}$ the equation for $Y$ is
\begin{equation*}
s + \sum_{i=1}^{3} (a_{i} s + b_{i}) y_{i}^{2}
\end{equation*}
where $y$ lies in the locus $s=0$ and the $a_{i}, b_{i}$ are constants satisfying:
\begin{itemize}
\item For any $i$ at most one of $a_{i}, b_{i}$ is zero.
\item Not all the $b_{i}$ vanish simultaneously.
\end{itemize}
The singularities will occur where
\begin{equation*}
(a_{i}s+b_{i})y_{i} = 0 \, \, \forall i \qquad \textrm{and} \qquad 1 + \sum a_{i}y_{i}^{2} = 0.
\end{equation*}
We show that a component of $\Sing(Y)$ containing $y$ has dimension $\leq 1$.  Note that it is impossible for none or all of the $b_{i}$ to vanish.  If WLOG only $b_{1}$ vanishes, then the component of $\Sing(Y)$ containing $y$ is defined by $s=y_{2}=y_{3} = 0$ and $1 + a_{1}y_{1}^{2} = 0$.  If WLOG only $b_{1}, b_{2}$ vanish then the singularities are defined by $s = y_{3} = 0$ and $1+y_{1}^{2} + y_{2}^{2} = 0$.  This verifies that $Y$ is normal.

Note that locally $Y$ will be defined by a cubic equation in $\mathbb{A}^{4}$.  Furthermore the equation will always involve a term of degree $\leq 2$ no matter which coordinate system we use, so that it does not describe a cone.  Using the classification of normal cubic surfaces (\cite{Sch}) we see that a general hyperplane section of $Y$ through any fixed point has only canonical singularities.  Thus $Y$ also has only canonical singularities by \cite[Theorem 5.34]{KM98}.
\end{proof}

By essentially the same argument we have:

\begin{lemm}
Let $\ell \subset \mathbb{P}^{3}$ be a plane and let $Y$ denote the $\pi_{1}$-preimage of $\ell$ in $X$.  Suppose that $Y \not \subset S$.  Then $Y$ is normal with canonical singularities.
\end{lemm}

\begin{proof}
In this case it is straightforward to check the claim about normality using a dimension count.  Fix an exceptional divisor $E$ over $Y$ with center $T$ of codimension $\geq 2$. Let $H \in |h_{1}|$ be a general element intersecting $T$.  Note that $K_{Y} + H$ is Cartier since it is the restriction of a Cartier divisor from the ambient variety.  Furthermore by Lemma \ref{lemm: threefoldquadricbundle} we know that $H$ is normal and has canonical singularities.  Applying inversion of adjunction (\cite[Corollary 1.4.5]{BCHM}), we see that the discrepancy of $(Y,H)$ with respect to $E$ is $\geq 0$.  Thus the discrepancy of $K_{Y}$ with respect to $E$ is also non-negative.
\end{proof}

\begin{coro} \label{coro: deg1section}
Let $\ell \subset \mathbb{P}^{3}$ be a linear subspace of dimension $1 \leq \dim(\ell) \leq 2$ and let $Y$ denote the $\pi_{1}$-preimage of $\ell$ in $X$.  Suppose that $Y \not \subset S$.  Then $a(Y,L) = 1$ and $(Y,L)$ is not adjoint rigid.
\end{coro}

\begin{proof}
Recall that $Y$ is normal, so by adjunction we have $K_{Y} + L|_{Y} = (3-\dim \ell)h_{1}|_{Y}$.  Since $Y$ has only canonical singularities this suffices to conclude the statement.
\end{proof}

\begin{prop} \label{prop: aclassificationforbhb}
If $Y \subset X$ is an adjoint rigid subvariety satisfying $a(Y,L) \geq 1$, then either:
\begin{enumerate}
\item $Y \subset S$, or
\item $Y$ is a smooth fiber of $\pi_{1}$, or
\item $Y$ is a line in a smooth fiber of $\pi_{1}$, or
\item $Y$ is a smooth fiber of $\pi_{2}$.
\end{enumerate}
Furthermore, every component $S' \subset S$ has $a(S',L) = 1$, we have $b(\overline{\mathbb{Q}},Y,L) = 2$ in case (2), and we have $b(\overline{\mathbb{Q}},Y,L) = 1$ in cases (3) and (4).
\end{prop}

\begin{proof}
We split into cases based on the dimension of $Y$.

\textbf{Dimension 4:} Let $Y$ be the resolution of an irreducible reduced divisor $Y'$ in $X$.  There are essentially three different options for the behavior of $Y$ under the two projection maps: either it can dominate both $\mathbb{P}^{3}$s, or it can map dominantly in one direction and onto a surface in the other.  We let $f_{1}$ and $f_{2}$ denote the composition of the map $Y \to Y'$ with the projection maps.  (Note that $f_{1}$ and $f_{2}$ need not have connected fibers.)

First suppose that $f_{1}: Y \to \mathbb{P}^{3}$ is dominant.  We will prove the claim by computing $H^{0}(Y,K_{Y} + 3h_{1} + 2h_{2})$.  Note that $2h_{1} + 2h_{2}$ is big and nef, so by Kawamata-Viehweg vanishing we have
\begin{equation*}
H^{0}(Y,K_{Y} + 3h_{1} + 2h_{2}) \geq H^{0}(Y_{1}, K_{Y_{1}} + 2h_{1} + 2h_{2})
\end{equation*}
where $Y_{1}$ is a general element of $|h_{1}|$.  Repeating this argument inductively several times, we see that
\begin{equation*}
H^{0}(Y,K_{Y} + 3h_{1} + 2h_{2}) \geq H^{0}(Y_{3}, K_{Y_{3}} + h_{1} + h_{2})
\end{equation*}
where $Y_{3}$ is the curve which is a complete intersection of $Y$ with two elements of $|h_{1}|$ and one element of $|h_{2}|$.  Since $Y_{3}$ is not contracted by either projection we have
\begin{equation*}
H^{0}(Y_{3}, K_{Y_{3}} + h_{1} + h_{2}) \geq 1.
\end{equation*}
This shows that $a(Y,L) \leq 1$.  In fact, this space of sections space will have dimension $\geq 2$ unless both $h_{1}$ and $h_{2}$ have degree $\leq 1$ on $Y_{3}$.  This means that $Y \cdot h_{1}^{3} \cdot h_{2} \leq 1$ and $Y \cdot h_{1}^{2} \cdot h_{2}^{2} \leq 1$, an impossibility.  So if $a(Y,L)=1$ then furthermore $Y$ can not be adjoint rigid.

The other case is when $Y$ maps to a surface under $f_{1}$ so that the class of $Y'$ is $ch_{1}$ for a positive integer $c$.  In this case we intersect $Y$ with a general member of $|h_{1}|$ and a general member of $|h_{2}|$ to obtain a smooth surface $Y_{2}$ satisfying
\begin{equation*}
H^{0}(Y,K_{Y} + 3h_{1} + 2h_{2}) \geq H^{0}(Y_{2}, K_{Y_{2}} + 2h_{1} + h_{2}).
\end{equation*}
Note that $h_{2}$ is a big and nef divisor on $Y_{2}$, so we may intersect against $|h_{1}|$ one more time while maintaining a surjection of sections.  This intersection is a (reducible) curve $Y_{3}$ which is a union of fibers of $f_{1}$.  Letting $C_{i}$ denote the fibers, we have
\begin{equation*}
H^{0}(Y,K_{Y} + 3h_{1} + 2h_{2}) \geq \sum_{i} H^{0}(C_{i}, K_{C_{i}} + h_{2}).
\end{equation*}
The right hand side can only vanish if $h_{2} \cdot C_{i} = 1$ for every $C_{i}$.  In this case, if $Y_{2} \to Z \to \mathbb{P}^{3}$ denotes the Stein factorization of $f_{1}$ then every fiber of $Y_{2} \to Z$ is mapped to a line in $\mathbb{P}^{3}$.  This means that for every fiber of the first projection map on $Y'$ the intersection against a general element of $h_{2}$ must be a union of lines.  But this is an impossibility, since the locus of fibers which are quadrics of corank $\geq 2$ are only parametrized by a one-dimensional space.  This shows that $a(Y,L) \leq a(X,L)$.  In particular we see that every component $S'$ of $S$ has $a(S',L) \leq 1$, and since $S'$ is covered by quadric cones it is clear the equality is achieved.

Next we show that $Y$ is not adjoint rigid.  Indeed, we can only have $\sum_{i} H^{0}(C_{i}, K_{C_{i}} + h_{2}) = 1$ if $Y_{3}$ is irreducible and has exactly degree $2$ against $h_{2}$.  In other words, we need $Y \cdot h_{1}^{2} \cdot h_{2}^{2} = 2$, showing that $c=1$.   This case is handled by Corollary \ref{coro: deg1section}.

\textbf{Dimension 3:} 
Let $Y$ be the resolution of an irreducible reduced threefold $Y'$ in $X$.  Just as before, we need to understand the possible behavior of $Y$ under the two projection maps.  In particular, we need to rule out the following two cases.  First, since both projection maps are flat on $X$ it is impossible for $Y$ to map to a point under either projection map.  Second, it is impossible for $Y$ to map to a curve $C$ under one projection map and to a surface $S$ under the other.  Indeed, in this case $Y'$ must be equal to $C \times S$, but no such subvariety is contained in $X$.

Just as in dimension $4$, we argue by computing $H^{0}(Y,K_{Y} + 3h_{1} + 2h_{2})$ by intersecting against general members of $|h_{1}|$ and $|h_{2}|$.  The only case that requires special discussion is when $Y$ maps to a curve under the first projection map and to $\mathbb{P}^{3}$ under the second projection map.  In this case $Y'$ must have numerical class $ch_{1}^{2}$ for some positive integer $c$.  Note that this is very similar to the second case in dimension $4$; by the same argument, we see that $a(Y,L) \leq 1$ unless possibly $Y \subset S$.  Furthermore, we see that either $a(Y,L) < 1$ or $Y$ is not adjoint rigid unless possibly $c=1$.  This final case is handled by Corollary \ref{coro: deg1section}.

\textbf{Dimension 2:}
The fibers $F_{2}$ of the second projection map are isomorphic to $\mathbb{P}^{2}$ and they each satisfy $a(F_{2},L) = a(X,L)$ and are adjoint rigid.  In this case $b(F_{2},L) = 1$.  The fibers $F_{1}$ of the first projection map are (possibly singular) quadrics.  When $F_{1}$ is smooth we have $a(F_{1},L)= 1$, $b(F_{2},L) = 2$ and the fiber is adjoint rigid.  When $F_{1}$ has corank $1$, $a(F_{1},L) = 1$, $b(F_{1},L) = 2$ and the fiber is adjoint rigid.  When $F_{1}$ has corank $2$, it is a union of two planes $P$.  Each plane satisfies $a(P,L) = 3/2$ and is adjoint rigid.  When $F_{1}$ has corank $3$, it is a double plane and the underlying reduced plane $P$ satisfies $a(P,L) = 3/2$ and is adjoint rigid.

Suppose now that $Y$ is the resolution of an irreducible reduced surface $Y'$ in $X$ and that $Y$ is not mapped to a point under either projection map.  The volume of $L$ is $9 Y \cdot h_{1}^{2} + 12 Y \cdot h_{1} \cdot h_{2} + 4 Y \cdot h_{2}^{2}$, and this is $>9$ under our assumption on $Y$.  So either $Y$ is not adjoint rigid or $a(Y,L)<1$ by Theorem \ref{theo: surfacebabbound}.

\textbf{Dimension 1:}
For any rational curve $C$ on $X$ we have $(3h_{1} + 2h_{2}) \cdot C \geq 2$ with equality if and only if $h_{1} \cdot C = 0$ and $h_{2} \cdot C = 1$.  So we have $a(C,L|_{C}) = a(X,L)$ precisely when $C$ is a line in a fiber of the first projection map.  In this case $C$ is adjoint rigid and $b(C,L|_{C}) = 1$.
\end{proof}

\subsection{Construction of the thin set}

We follow Steps (0)-(9) of Section \ref{sect: computation}.  For Step (0), note that $\mathbf{B}_{+}(L)$ is empty.

Step (1) is more-or-less accomplished by Proposition \ref{prop: aclassificationforbhb}.  Precisely speaking, there is still a small ambiguity about subvarieties of $S$.  However, we will remove all rational points contained in $S$ anyways, so this minor issue will be irrelevant for the final computation.  

In Step (2) we make no changes.  In Step (3) we set $V=S$, and we will remove the rational points on $V$ later on.  In Step (4) we need to consider the universal family of lines in fibers of the first projection.  Note that the family of lines has $b$-value smaller than the $b$-value of $X$.  Furthermore, there is no adjoint rigid cover of $\mathbb{P}^{1}$ with the same $a$-value.  Thus there is no need to remove any points from this family.

For Steps (5)-(8), we must consider the fibers of the two projection maps and adjoint rigid covers of $X$.  Consider the second projection.  Here the fibers are isomorphic to $\mathbb{P}^{2}$; in particular $b(\overline{\mathbb{Q}},Y,L)$ is only $1$.  Furthermore, there are no covers of projective space with the same $a$-value.  Thus we do not need to remove any points coming from this family.  

Now consider the first projection whose fibers are quadrics $Q$.  The first step is to shrink the family so that it is smooth; thus we increase $V$ by adding all singular fibers of $\pi_{1}$.  Let $W^{\circ} \subset \mathbb{P}^{3}$ denote the complement of the coordinate hyperplanes.  We then take a base change via the map $\mathbb{P}^{3} \to \mathbb{P}^{3}$ defined in coordinates by $(x_0:x_1:x_2:x_3) \mapsto (x_0^2: -x_1^2: x_2^2:-x_3^2)$.  Let $W'$ denote the preimage of $W^{\circ}$.  This map is Galois, kills the geometric monodromy on the N\'eron-Severi spaces of the fibers, and the family induced by base change admits a section over the ground field.  Note that the map $W' \to W^{\circ}$ is \'etale.  In step (6), there are no non-trivial adjoint rigid covers of $Q$ with the same $a$-value, so we only need to consider the identity map.  Thus Step (7) yields a unique family, namely the original family $p: \mathcal{U}' \to W'$.  In Step (8) there is no need to take a further base-change.  Call the closure of the resulting family $q: \mathcal{Y} \to T$ with maps $f: \mathcal{Y} \to X$ and $g: T \to \mathbb{P}^{3}$.  Now the $b$-value of $\mathcal{Y}$ is $2$, so we must remove $f(\mathcal{Y}(\mathbb{Q}))$.  This will exactly consist of the points on the fibers over $g(T(\mathbb{Q}))$; note that all of these fibers are split.  As we vary over twists of $f: \mathcal{Y} \to X$, we remove the points only from the twists with $b$-value equal to $2$, and their union will consist of the points on the split fibers of $\pi_{1}$.  We also must remove contributions from the branch locus, but this is already contained in $V$.

Finally, we must show that there are no dominant thin maps $f: Y \to X$ such that $a(Y,f^{*}L) = a(X,L)$ and $(Y,f^{*}L)$ is adjoint rigid.  By Theorem \ref{theo: akash} it suffices to show that every divisor $D$ on $X$ has $a(D,L) \leq a(X,L)$.  This was proved in Proposition \ref{prop: aclassificationforbhb}.

Note that a diagonal quadric surface $\sum a_{i} y_{i}^{2}$ that carries a rational point will split if and only if $a_{0}a_{1}a_{2}a_{3}$ is a square. 
Altogether we find:

\begin{theo}
The conjectural exceptional set $Z$ is the union of points on fibers $X_{p}$ of the first projection such that $p = (a_{0}:a_{1}:a_{2}:a_{3}) \in \mathbb{P}^{3}$ and either $a_{0}a_{1}a_{2}a_{3}$ is a non-zero square or some $a_{i}=0$.
\end{theo}

\cite[Theorem 1.1]{BHB18} proves that Manin's Conjecture holds with Peyre's constant after removing precisely this exceptional set.

\section{Open questions} \label{sect: openques}

We finish with a list of open questions.

\subsection{Behavior of a and b invariants}

\begin{ques}
Is the $a$-invariant constant in smooth uniruled families polarized by a big and nef divisor?  (See Theorem \ref{theo: aconstantinfamilies}.)
\end{ques}

\begin{ques}
Suppose our ground field is algebraically closed.  Is the $b$-constant lower semi-continuous in smooth uniruled families? (See Theorem \ref{theo: bconstantinfamilies}.)  
\end{ques}

\begin{ques} \label{ques: etalecodim1}
Let $f: Y \to X$ be a dominant generically finite map of smooth projective varieties and let $L$ be a big and nef divisor on $X$ such that $a(Y,f^{*}L) = a(X,L)$ and $(Y,f^{*}L)$ is adjoint rigid.  Is it true that $f$ is birationally equivalent to a map of klt varieties that is \'etale in codimension $1$?
\end{ques}

\cite{LTDuke} conjectures an affirmative answer to Question \ref{ques: etalecodim1}.  If the question indeed has an affirmative answer, then one can identify all such covers $Y$ by computing the \'etale fundamental group of a suitably chosen singular birational model $X'$ of $X$.  Using the finiteness of the \'etale fundamental group of $\mathbb{Q}$-Fanos proved by \cite{Xu14}, one would obtain an alternative approach to the classification problem of adjoint rigid covers with the same $a$-value.  The advantage of this viewpoint is that it is somewhat more explicit, allowing one to construct the covers directly rather than going through the round-about procedure in Step (6) of Section \ref{sect: computation}.

\subsection{Examples}

In the text we gave several classes of examples where we could classify all thin maps $f: Y \to X$ with larger $a$ and $b$ invariants.  Here are several examples we do not know how to handle.

\begin{ques}
Let $X$ be a smooth hypersurface in $\mathbb{P}^{n}$ of degree $\leq n/2$ and let $L$ denote the hyperplane class.  Is it true that there is no thin map $f: Y \to X$ such that $a(Y,f^{*}L) = a(X,L)$?
\end{ques}

For surfaces, the following two questions are the key for understanding the thin exceptional set.

\begin{ques}[= Question \ref{ques: ratcomponents}]
Let $X$ be a geometrically integral smooth del Pezzo surface of degree $1$ over a number field.  Let $M$ denote the curve parametrizing rational curves in $|-2K_{X}|$.  Does every component of $M$ have genus $\geq 2$?
\end{ques}

\begin{ques}[= Question \ref{ques: weakdelpezzo}]
Let $X$ be a geometrically integral smooth weak del Pezzo surface over a number field and let $L = -K_{X}$.  Suppose that $f: Y \to X$ is a dominant generically finite map such that $a(Y,f^{*}L) = a(X,L)$ and $(Y,f^{*}L)$ is adjoint rigid.  Then is $b(F,Y,f^{*}L) < b(F,X,L)$?
\end{ques}

Note that the analogue of Question \ref{ques: weakdelpezzo} over an algebraically closed field is solved by \cite[Section 6]{LTDuke}, so the difficulty arises from working over a number field.

\begin{ques}
Let $X$ be a toric variety.  Give a combinatorial description of all thin maps $f: Y \to X$ such that $a(Y,f^{*}L) \geq a(X,L)$.
\end{ques}

See Example \ref{exam: toricvarieties} for more information on the $a$ and $b$ invariants for toric varieties.

\begin{ques}
Let $X \subset \mathbb{P}^{n} \times \mathbb{P}^{m}$ be a hypersurface of bidegree $(1,2)$ or $(1,3)$.  Can one give a uniform analysis of the $a,b$-invariants for some class of such varieties?
\end{ques}

These hypersurfaces should provide a rich source of examples where the exceptional set is Zariski dense (as in \cite{BT-cubic} and \cite{BHB18}).

\nocite{*}
\bibliographystyle{alpha}
\bibliography{absurvey}

\end{document}